\documentclass[final,1p,times]{elsarticle}

\usepackage{amssymb}
\usepackage{amsmath}
\usepackage{hyperref}

\journal{Journal of Computational Physics}

\renewcommand{\aa}{\mathbf{a}}

\newcommand{\cc}{\mathbf{c}}
\newcommand{\dd}{\mathbf{d}}
\newcommand{\ee}{\mathbf{e}}
\newcommand{\ff}{\mathbf{f}}
\newcommand{\kk}{\mathbf{k}}
\newcommand{\nn}{\mathbf{n}}
\newcommand{\pp}{\mathbf{p}}

\newcommand{\rr}{\mathbf{r}}
\renewcommand{\ss}{\mathbf{s}}
\renewcommand{\tt}{\mathbf{t}}

\newcommand{\xx}{\mathbf{x}}
\newcommand{\yy}{\mathbf{y}}
\newcommand{\zz}{\mathbf{z}}

\renewcommand{\AA}{\mathbf{A}}
\newcommand{\BB}{\mathbf{B}}

\newcommand{\EE}{\mathbf{E}}
\newcommand{\FF}{\mathbf{F}}
\newcommand{\GG}{\mathbf{G}}
\newcommand{\HH}{\mathbf{H}}
\newcommand{\II}{\mathbf{I}}

\newcommand{\PP}{\mathbf{P}}
\newcommand{\QQ}{\mathbf{Q}}
\newcommand{\VV}{\mathbf{V}}
\newcommand{\WW}{\mathbf{W}}
\newcommand{\ZZ}{\mathbf{Z}}
\newcommand{\bzero}{\mathbf{0}}
\newcommand{\btau}{\boldsymbol\tau}
\newcommand{\bsigma}{\boldsymbol\sigma}
\newcommand{\brho}{\boldsymbol\rho}
\newcommand{\bxi}{\boldsymbol\xi}

\newcommand{\cW}{\mathcal{W}} 
\newcommand{\cM}{\mathcal{M}}
\newcommand{\cN}{\mathcal{N}}

\newcommand{\cP}{\mathcal{P}}
\newcommand{\cL}{\mathcal{L}}
\newcommand{\cR}{\mathcal{R}}

\newcommand{\bbC}{\mathbb{C}} 
\newcommand{\bbR}{\mathbb{R}} 
\newcommand{\bbZ}{\mathbb{Z}} 

\newcommand{\iu}{\mathrm{i}} 

\newcommand{\inc}{\mathrm{inc}}

\renewcommand{\paragraph}[1]{\vspace{0.1in}\noindent\textbf{#1}}

\begin{document}

\begin{frontmatter}



\title{Method of Fundamental Solutions for Maxwell's Equations in Bi-Periodic Multilayered Media}

 \author[label1]{Jared Weed}
 \author[label1]{Bowei Wu}
 \author[label2]{Jingfang Huang}
 \author[label1]{Min Hyung Cho\corref{cor1}}
 \cortext[cor1]{Corresponding Author}
 \ead{minhyung_cho@uml.edu}
 \affiliation[label1]{organization={Department of Mathematics and Statistics, University of Massachusetts Lowell},
             city={Lowell},
             postcode={01854},
             state={MA},
             country={USA}}

 \affiliation[label2]{organization={Department of Mathematics, University of North Carolina at Chapel Hill},
             city={Chapel Hill},
             postcode={27599},
             state={NC},
             country={USA}}


\begin{abstract}
In this paper, we present an accurate numerical method for the time-harmonic Maxwell's equations for bi-periodic multilayered media with quasi-periodic incident waves using the Method of Fundamental Solutions in conjunction with a periodization scheme. Following an approach used in acoustic scattering problems, the electric and magnetic fields in each layer are expressed as a sum of near and distant interactions. The near interaction comprises interactions between the unit cell and its nearest neighboring copies, while the distant interaction is approximated by proxy source points placed on spheres surrounding the unit cell. Imposing continuity of tangential components at the layer interface, quasi-periodicity conditions on the walls of the unit cell, and Rayleigh-Bloch expansion for the radiation condition yields a system of equations for the unknown coefficients, which can be solved by Schur complement and a backward-stable solver. The scheme is verified with known solutions and exhibits exponential convergence close to $10^{-14}$ for both single and multiple interfaces. An example with 39 interfaces is presented to demonstrate the solver's performance. The paper provides promising results for extending this method to a fast and accurate boundary integral equation solver for many cutting-edge applications involving a large number of layers in electromagnetics and optics.
\end{abstract}




\end{frontmatter}



\section{Introduction}\label{sec:intro}
Electromagnetic (EM) waves in multilayered media play crucial roles in modern high-tech devices such as semiconductor packaging \cite{deboi2022improved,kim2022emi}, diffraction gratings, thin-film photovoltaics \cite{atwater2011plasmonics,kelzenberg2010enhanced}, and high-power lasers \cite{perry1995high,barty2004overview}. As these devices become more complex, accurate numerical simulations of electromagnetic fields in multilayered media require efficient algorithms and significant computational resources to optimize performance. Many commercial and in-house software packages for electromagnetic simulations employ methods such as the finite element method \cite{monk2003finite,bao1995finite}, finite difference method \cite{taflove2005computational,yee1966numerical}, method of moments \cite{harrington1993field}, rigorous coupled-wave analysis \cite{moharam1981rigorous,li1996use,cho2008rigorous}, boundary integral equation (BIE) method \cite{bruno2016windowed,bruno2017windowed,PerezArancibia2019DomainDF,nicholls2020sweeping,cho2012parallel,chen2018accurate,lai2014fast,cho2015robust,zhang2021fast,cho2019spectrally,boag1992analysis,cho2021adapting,wang2026high}, and method of fundamental solutions (MFS) \cite{kupradze1964method}. Each of these approaches has advantages and disadvantages, and an exhaustive comparison is beyond the scope of this paper.

Green’s function-based methods, such as the BIE method and the MFS, offer several well-known benefits for exterior domain problems, including automatic enforcement of radiation conditions and dimensionality reduction. Traditionally, BIE methods represent the scattered field using layer potentials with unknown density functions defined on the boundary. Enforcing boundary conditions yields integral equations for these unknown densities, which are then solved numerically using appropriate quadrature rules. Accurate quadrature is a critical component for the success of integral equation methods and remains an active research area \cite{wu2021zeta,weed2023quadrature,ding2021quadrature}.

The MFS operates similarly: the scattered field is represented as a linear combination of Green's functions evaluated at target points influenced by artificial source points offset from the boundary. The unknown coefficients are determined by enforcing boundary conditions on this representation. Because source points are offset, singularities in the Green's function are avoided, making the MFS simpler to implement than BIE methods. However, careful placement of source points is essential to achieve accurate solutions.

In this paper, we consider EM scattering in multilayered media where the interfaces between layers are bi-periodic and the incident wave is quasi-periodic. There are two main approaches to apply the MFS or BIE methods in multilayered media: (a) the layered media Green’s function method and (b) the free-space Green’s function method with a periodization scheme. In the first approach, the layered media (or background) Green’s function, which satisfies the layer interface conditions, must be determined, after which the scattered waves are computed as if in free space. However, accurately computing the layered media Green’s function is often challenging and computationally expensive \cite{cho2021adapting,cho2012parallel,cho2017efficient,cho2018heterogeneous,chen2018accurate}.

In the second approach, the solution in each layer can be expressed using a quasi-periodic Green’s function, which is known to suffer from convergence issues. In this work, we employ a periodization scheme, previously shown to be effective for both the two- and three-dimensional Helmholtz equations \cite{cho2015robust}, to overcome these issues in the context of the time-harmonic Maxwell's equations. The resulting linear system, obtained by matching continuity conditions at the layer interfaces, can be solved efficiently by eliminating the additional unknowns introduced by the periodization scheme via the Schur complement. Similar techniques have been applied to acoustic wave scattering problems governed by the Helmholtz equation in both two \cite{cho2015robust} and three dimensions \cite{cho2019spectrally,wu2023robust,liu2016numerical}.

Section \ref{sec:scatform} describes the boundary value problem for a free-space EM scattering problem with the transmission boundary condition. Section \ref{sec:single} presents the formulation of the periodization scheme for a single interface. Section \ref{sec:multi} extends the single-interface formulation to general multilayered media with minimal additional effort. In Section \ref{sec:experiments}, we present numerical experiments across a range of parameters to examine the performance of the numerical method. Finally, Section \ref{sec:conclusion} concludes the paper, with a discussion of directions for future work.


\begin{figure}[t] 
   \centering
   \includegraphics[width=2in]{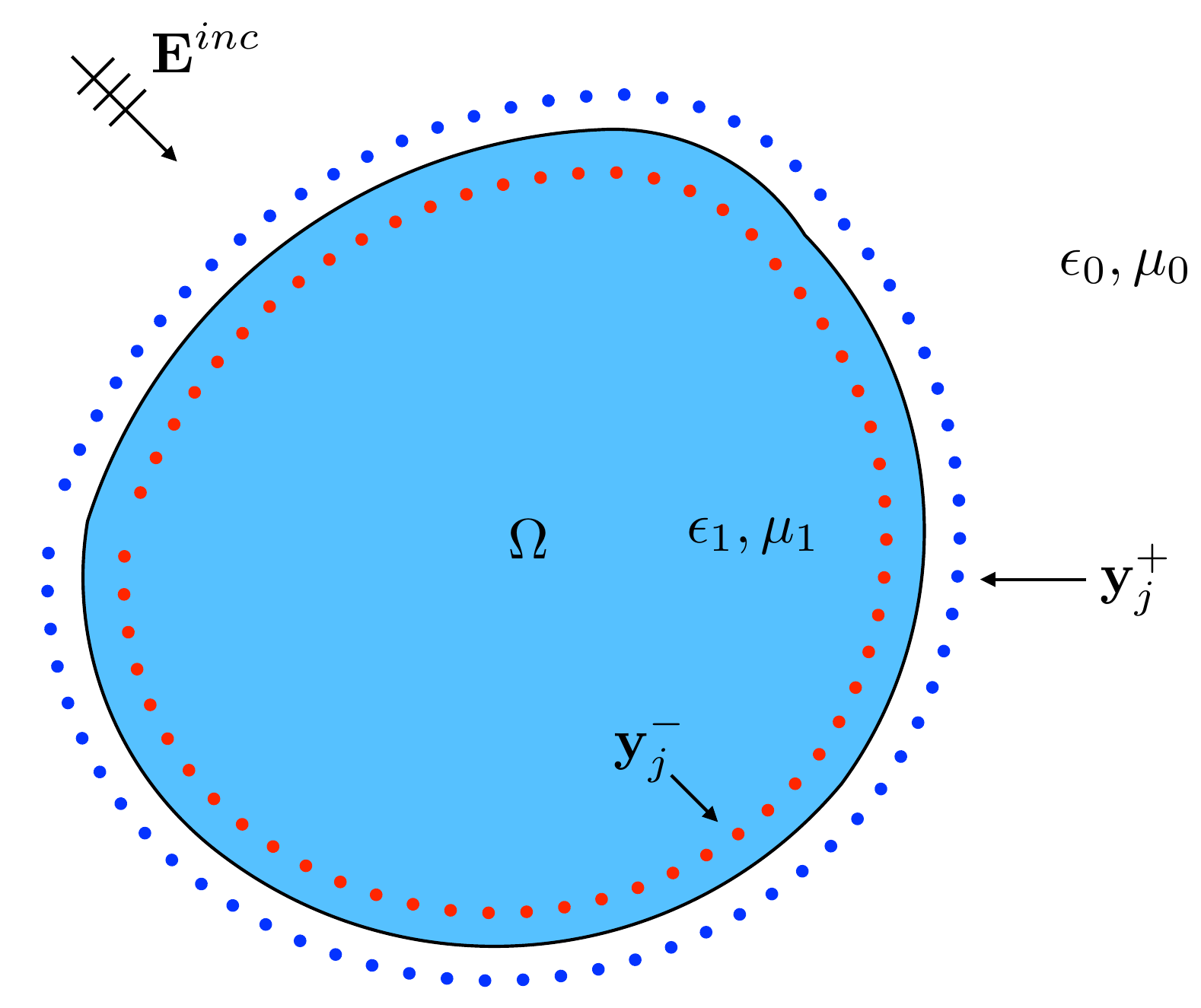} 
   \caption{Free Space Problem Setup}
   \label{fig:free}
\end{figure}
\section{Free-Space EM scattering MFS Formulation}\label{sec:scatform}
Consider a scattering object $\Omega$ with permittivity and permeability $\epsilon_1$ and $\mu_1$ embedded in free space with permittivity $\epsilon_0$ and $\mu_0$ (see Fig. \ref{fig:free}). The electric and magnetic fields outside ($\EE_0$, $\HH_0$) and inside ($\EE_1$, $\HH_1$) the object resulting from the incident wave $\EE^{\inc}(\xx) = \EE^{\inc}e^{\iu \kk \cdot \xx}$ with the wavevector $\kk = (k_x, k_y, k_z)$ scattered by $\Omega$ can be found by solving the time-harmonic Maxwell's equations
\begin{align}
    \nabla\times\EE_l &= \iu \omega \mu_l \mathbf{H}_l,\\
    \nabla\times\HH_l &= -\iu \omega \epsilon_l \mathbf{E}_l, \label{maxwell:2}\\
    \nabla \cdot\EE_l &=0,\\
    \nabla \cdot\HH_l &=0,
\end{align}
with continuity of tangential fields at the boundary of $\Omega$,
\begin{align}
\nn\times \EE_1 =  \nn\times \EE_0+\nn\times\EE^{\inc} \mbox{ and }~ \nn\times \HH_1 =  \nn\times \HH_0,\label{eq:bc2}
\end{align}
and a Silver-M\"uller radiation condition, where $l=0,1$ and $\mathbf{n}$ is the outward unit normal vector.  In the framework of MFS, the electric and magnetic fields inside ($l = 1$) and outside ($l = 0$) of $\Omega$ are represented by
\begin{align}
&    \EE_0(\xx) \approx \sum_{j=1}^\cN c^{+}_j \GG^0_\EE(\xx,\yy^{-}_j,\pp^{-}_j) \quad \mathrm{and} \quad \HH_0(\xx) \approx \sum_{j=1}^\cN c^{+}_j \GG^0_\HH(\xx, \yy^{-}_j, \pp^{-}_j), ~~ \textrm{for } \xx \in \mathcal{R}^{3}\backslash \Omega, \\
&    \EE_1(\xx) \approx \sum_{j=1}^\cN c^{-}_j \GG^1_\EE(\xx,\yy^{+}_j,\pp^+_j) \quad \mathrm{and} \quad \HH_1(\xx) \approx \sum_{j=1}^\cN c^{-}_j \GG^1_\HH(\xx, \yy^{+}_j, \pp^+_j), ~~ \textrm{for } \xx \in \Omega, \label{mfsoffset}
\end{align}
where $\mathcal{N}$ artificial dipole source points are placed inside ($\yy_j^{-}$) and outside ($\yy_j^{+}$) of the boundary of the domain oriented along $\pp^{-}_j$ and $\pp^{+}_j$, respectively. The electric and magnetic free-space Green's functions with wavenumber $k_l=\omega\sqrt{\epsilon_l \mu_l}$ are given by
\begin{equation}
    \GG^l_\EE(\xx,\yy,\pp) :=\frac{\iu\omega}{4\pi}\left(\II + \frac{\nabla\nabla}{k_l^2}\right)\frac{e^{\iu k_l r}}{r}\cdot\pp \quad \mathrm{and} \quad \GG^l_\HH(\xx,\yy,\pp) =\frac{1}{4\pi\mu_l}e^{\iu k_l r}\frac{\iu k_l r-1}{r^3}\rr\times\pp,   
\end{equation}
where $\rr= \xx-\yy$ and $r=||\xx-\yy||$. By applying the boundary conditions in Eq. (\ref{eq:bc2}) using two perpendicular tangential vectors on the boundary with source dipoles, the unknown coefficients $c^{\pm}_j$ can be determined. 

A similar idea can be applied to an EM transmission problem for a bi-periodic surface with periodicity $d_x$ and $d_y$ along the $x$- and $y$-directions by replacing the free-space Green's function with the quasi-periodic Green's function
\begin{equation}
\GG^{QP}_\FF(\xx, \yy, \pp) =   \sum_{m,n\in\bbZ}\alpha^m\beta^n\GG^l_\FF(\xx+m\ee_x, \yy+n\ee_y,\pp), ~~\textrm{for } \FF = \EE \mbox{ or } \HH,
\end{equation}
where $\alpha=e^{\iu k_x d_x}$ and $\beta=e^{\iu k_y d_y}$ are the Bloch phases along the $x$- and $y$-directions, and $\ee_x = (d_x, 0, 0)$ and $\ee_y = (0, d_y, 0)$.\footnote{A function $\ff:\bbR^3 \to\bbC^3$ is \textit{quasi-periodic} if $\ff(x,y,z)=\alpha^{-1}\ff(x+d_x,y,z)=\beta^{-1}\ff(x,y+d_y,z)$} The scattered field must satisfy the outward radiation condition through Rayleigh-Bloch expansions at a certain distance away from the surface, namely,
\begin{align}
\EE(\xx) &= \sum_{m,n \in \bbZ} \aa_u^{mn}\exp(\iu(\kappa_x^m x + \kappa_y^n y + \kappa_{0,z}^{mn}(z-z_u))), ~~\textrm{for } z \geq z_u, \label{eq:radation1}\\
\EE(\xx) &= \sum_{m,n \in \bbZ} \aa_d^{mn}\exp(\iu(\kappa_x^m x + \kappa_y^n y + \kappa_{1,z}^{mn}(-z+z_d))), ~~\textrm{for } z \leq z_d,\label{eq:radation2}
\end{align}
where $\kappa_x^m = k_x+2\pi m/ d_x$, $\kappa_y^n = k_y+2\pi n/ d_y$, and $\kappa_{l, z}^{mn} = +\sqrt{k_l^2-(\kappa_x^m)^2-(\kappa_y^n)^2}$, for $l=0,1$.

It is, however, well-known that the quasi-periodic Green's function suffers from slow convergence and other non-physical issues. Therefore, a periodization scheme used for the three-dimensional Helmholtz equation is employed by approximating the quasi-periodic Green's function by near and distant field \cite{cho2019spectrally,wu2023robust,liu2016numerical}.
\begin{figure}[htbp] 
   \centering
   \includegraphics[width=4.4in]{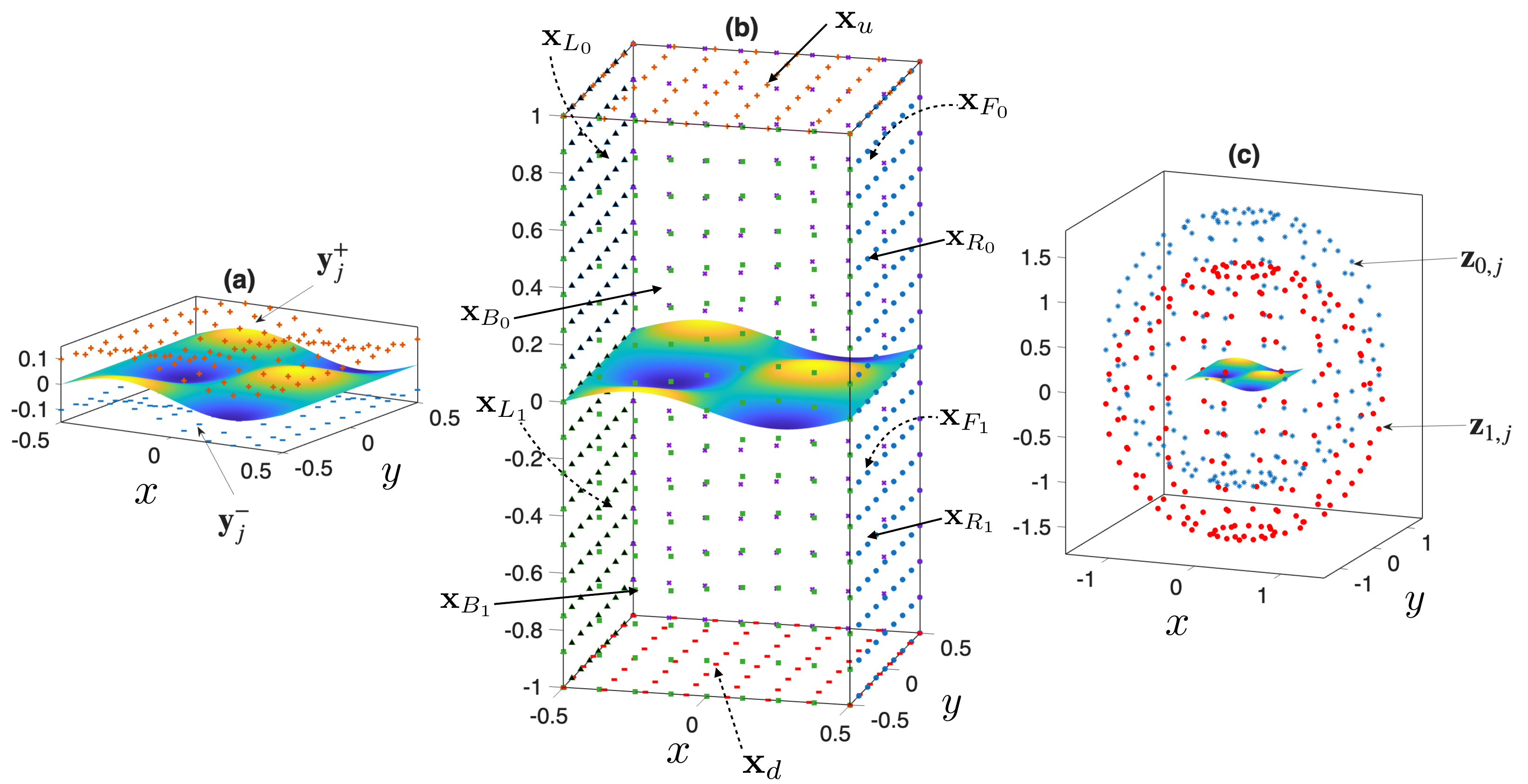} 
   \caption{(a) Unit cell and artificial source points above ($\yy_j^{+}$) and below ($\yy_j^-$) the layer interface $\Gamma$. (b) Left, right, back, and front walls in both top and bottom layers to enforce quasi-periodicity, and artificial boundary in the top $\xx_u$ and bottom $\xx_d$ layers to enforce the radiation condition. (c) Proxy source points for the top ($\zz_{0,j}$) and bottom ($\zz_{1,j}$) layers.}
   \label{fig:domain_one}
\end{figure}

\section{MFS for a single bi-periodic interface} \label{sec:single}
Consider an interface $\Gamma$ defined by a function $z=f(x,y)$ in the $x$-$y$ plane that is bi-periodic within the unit cell $[-\frac{d_x}{2},\frac{d_x}{2}]\times[-\frac{d_y}{2}, \frac{d_y}{2}]$---or more precisely, $f(x,y)=f(x+m d_x, y+n d_y) \; \forall m,n\in \bbZ$. Let $\Omega_0$ with relative permeability $\mu_0$ and relative permittivity $\epsilon_0$ and $\Omega_1$ with relative permeability $\mu_1$ and relative permittivity $\epsilon_1$ denote the domains above and below the interface $\Gamma$, respectively. We refer to $\Omega_0$ and $\Omega_1$ as the first and second layers, respectively; the definition extends naturally to multiple interfaces (see Section \ref{sec:multi}). Using the periodization idea from the Helmholtz equation, the scattered electric and magnetic fields in each domain can be expressed as a sum of near and distant interactions through proxy source points. 
The problem setting shown in Fig. \ref{fig:domain_one} is summarized as follows.

Let $\{\yy_j^{\pm}\}_{j=1}^\cN$ be MFS source points above ($+$) and below ($-$) the interface, and $\btau_j^{\pm}$ and $\bsigma_j^{\pm}$ are two tangential vectors at $\yy^{\pm}_j$. We choose $\yy_j^{\pm}$ by discretizing an $\sqrt{\cN}\times\sqrt{\cN}$ set of uniform collocation points on the interface and extending them by a factor of their unit normal in the positive and negative directions. For $l=0$ (top layer) and $l=1$ (bottom layer), let $\{\zz_{l,j}\}_{j=1}^\cP$ be $\cP$ proxy source points on a sufficiently-large spherical surface of radius $R$ centered on and enclosing the unit cell of $\Omega_l$, and $\brho_{l,j}$ and $\bxi_{l,j}$ are two tangential vectors at the proxy source points $\zz_{l,j}$. We choose $\zz_{l,j}$ by projecting a Fibonacci spiral of $\cP$ points in $[0,1)^2$ onto the surface of the unit sphere and extending each point by $R$. The unit cell of $\Omega_l$ is surrounded by left, right, back, and front artificial side walls denoted by $\xx_{L_l}$, $\xx_{R_l}$, $\xx_{B_l}$, and $\xx_{F_l}$, respectively. Finally, the first and second layers are also enclosed by artificial walls $\xx_u$ (above) and $\xx_d$ (below), respectively, on which the radiation condition will be imposed. Based on these notations, the scattered electric fields in the first and the second layers can be represented by
\begin{align}
    \EE_0(\xx) &\approx  \sum_{j=1}^{\cN} c^{\btau-}_{j} \tilde{\GG}^0_\EE(\xx,\yy^{-}_{j},\btau^-_j) + \sum_{j=1}^{\cN} c^{\bsigma-}_{j} \tilde{\GG}^0_\EE(\xx,\yy^{-}_{j},\bsigma^-_j)\nonumber\\
    &~~~~~~~~~~~~+\sum_{j=1}^{\cP}d^{\brho}_{0,j}\GG^0_\EE(\xx,\zz_{0,j},\brho_{0,j})+\sum_{j=1}^{\cP}d^{\bxi}_{0,j}\GG^0_\EE(\xx,\zz_{0,j},\bxi_{0,j}),~~\textrm{for } \xx \in \Omega_0 ,\label{eq:periodic1} \\
    \EE_1(\xx) &\approx  \sum_{j=1}^{\cN} c^{\btau+}_{j} \tilde{\GG}^1_\EE(\xx,\yy^{+}_{j},\btau^+_j) + \sum_{j=1}^{\cN} c^{\bsigma+}_{j} \tilde{\GG}^1_\EE(\xx,\yy^{+}_{j},\bsigma^+_j) \nonumber\\
    &~~~~~~~~~~~~+\sum_{j=1}^{\cP}d^{\brho}_{1,j}\GG^1_\EE(\xx,\zz_{1,j},\brho_{1,j})+\sum_{j=1}^{\cP}d^{\bxi}_{1,j}\GG^1_\EE(\xx,\zz_{1,j},\bxi_{1,j}),~~\textrm{for } \xx \in \Omega_1 ,\label{eq:periodic2} 
 \end{align}
where
\begin{align}
\tilde{\GG}^l_\EE(\xx,\yy,\pp) &= \sum_{\substack{m,n\in\bbZ \\ |n|,|m|\leq 1}}\alpha^m\beta^n\GG^l_\EE(\xx+m\ee_x, \yy+n\ee_y,\pp),~~\textrm{for } l=0,1.
\end{align}
By symmetry, the magnetic field in each layer can be expressed by replacing $\EE$ with $\HH$ and $\GG_\EE$ with $\GG_\HH$. Therefore, the magnetic field formulas have been omitted. The first two summations within $\EE_{0,1}(\xx)$ represent the near interaction between the unit cell and its immediate eight neighboring copies while the last  two  approximate the interaction between the unit cell and distant copies. Below, the continuity of tangential components of the electric and magnetic fields at the layer interface, quasi-periodicity condition on the sides of the unit cell, and the radiation condition at artificial walls above and below the layer interface are enforced and are represented in matrix form.

\paragraph{Transmission condition:} Let $\xx$ be a target point on the layer interface $\Gamma$ with $\tt$ and $\ss$ denoting two tangential vectors at $\xx$. The continuity of tangential components of the electric and magnetic fields at the interface can be expressed as
\begin{align}
\tt \cdot \EE_0(\xx) - \tt \cdot \EE_1(\xx) &= -\tt \cdot \EE^{\inc}(\xx),\\
\ss \cdot \EE_0(\xx) - \ss \cdot \EE_1(\xx) &= -\ss \cdot  \EE^{\inc}(\xx),\\
\tt \cdot \HH_0(\xx) - \tt \cdot \HH_1(\xx) &= \mathbf{0},\\
\ss \cdot \HH_0(\xx) - \ss \cdot \HH_1(\xx) &= \mathbf{0}.
\end{align}
By evaluating the electric field equations (\ref{eq:periodic1}) and (\ref{eq:periodic2}) and the corresponding magnetic field formulas at $\xx$ and substituting them into the boundary conditions, we obtain a system of equations for unknown coefficient vectors $\cc^{\btau\pm} = \{c^{\btau\pm}_{j}\}_{j=1}^\cN$, $\cc^{\bsigma\pm} = \{c^{\bsigma\pm}_{j}\}_{j=1}^\cN$, $\dd^{\brho}_0 = \{d^{\brho}_{0,j}\}_{j=1}^\cP$, $\dd^{\brho}_1 = \{d^{\brho}_{1,j}\}_{j=1}^\cP$, $\dd^{\bxi}_0 = \{d^{\bxi}_{0,j}\}_{j=1}^\cP$, and $\dd^{\bxi}_1 = \{d^{\bxi}_{1,j}\}_{j=1}^\cP$. For the sake of simplicity, only the first equation $\tt \cdot \EE_0(\xx) - \tt \cdot \EE_1(\xx) = -\tt \cdot \EE^{\inc}(\xx)$ is expanded here as
\begin{align}
&\tt \cdot \sum_{j=1}^{\cN} c^{\btau-}_{j} \tilde{\GG}^0_\EE(\xx,\yy^{-}_{j},\btau^-_j) + \tt \cdot  \sum_{j=1}^{\cN} c^{\bsigma-}_{j} \tilde{\GG}^0_\EE(\xx,\yy^{-}_{j},\bsigma^-_j)\nonumber\\
-&\tt \cdot \sum_{j=1}^{\cN} c^{\btau+}_{j} \tilde{\GG}^1_\EE(\xx,\yy^{+}_{j},\btau^+_j) -\tt \cdot  \sum_{j=1}^{\cN} c^{\bsigma+}_{j} \tilde{\GG}^1_\EE(\xx,\yy^{+}_{j},\bsigma^+_j)\nonumber \\
+&\tt \cdot  \sum_{j=1}^{\cP}d^{\brho}_{0,j}\GG^0_\EE(\xx,\zz_{0,j},\brho_{0,j})+\tt \cdot \sum_{j=1}^{\cP}d^{\bxi}_{0,j}\GG^0_\EE(\xx,\zz_{0,j},\bxi_{0,j})\nonumber\\
-&\tt\cdot \sum_{j=1}^{\cP}d^{\brho}_{1,j}\GG^1_\EE(\xx,\zz_{1,j},\brho_{1,j})-\tt\cdot \sum_{j=1}^{\cP}d^{\bxi}_{1,j}\GG^1_\EE(\xx,\zz_{1,j},\bxi_{1,j}) = -\tt_i\cdot \EE^{\inc}(\xx).
\end{align}
The second equation can be obtained by the substitution of $\tt$ with $\ss$. The third and fourth equations can be obtained by substituting $\EE$ with $\HH$ in the first and second equations and equating them with $\mathbf{0}$. For a set of $\cM$ target points $\{\xx_i\} _{i=1}^{\cM} \in \Gamma$ with tangent vectors $\{\tt_i\}_{i=1}^{\cM}$ and $\{\ss_i\}_{i=1}^\cM$, all four transmission boundary conditions can be represented in a matrix form equation
$$\left[\begin{array}{cc} \AA_{0,0}  & \AA_{0,1} \end{array}\right]\left[\begin{array}{c}  \cc^- \\  \cc^+ \\ \end{array}\right]+\left[\begin{array}{cc} \BB_{0,0}  & \BB_{0,1} \end{array}\right]\left[\begin{array}{c}  \dd_0 \\  \dd_1\\ \end{array}\right] = \ff.$$
For notational simplicity, we define
\begin{equation}
\tilde{\FF}^l(\xx, \yy, \btau, \bsigma) := \left[\begin{array}{cc}
\tt \cdot  \tilde{\GG}^l_\EE(\xx,\yy,\btau)  &  \tt \cdot  \tilde{\GG}^l_\EE(\xx,\yy,\bsigma) \\
\ss \cdot  \tilde{\GG}^l_\EE(\xx,\yy,\btau) &  \ss \cdot  \tilde{\GG}^l_\EE(\xx,\yy,\bsigma)\\
\tt \cdot  \tilde{\GG}^l_\HH(\xx,\yy,\btau)  & \tt \cdot  \tilde{\GG}^l_\HH(\xx,\yy,\bsigma)  \\
\ss \cdot  \tilde{\GG}^l_\HH(\xx,\yy,\btau) &  \ss \cdot  \tilde{\GG}^l_\HH(\xx,\yy,\bsigma)\end{array}\right],
\end{equation}
and
\begin{equation}
\FF^l(\xx, \yy, \btau, \bsigma) := \left[\begin{array}{cc}
\tt \cdot  {\GG}^l_\EE(\xx,\yy,\btau)  &  \tt \cdot  {\GG}^l_\EE(\xx,\yy,\bsigma) \\
\ss \cdot  {\GG}^l_\EE(\xx,\yy,\btau) &  \ss \cdot  {\GG}^l_\EE(\xx,\yy,\bsigma)\\
\tt \cdot  {\GG}^l_\HH(\xx,\yy,\btau)  & \tt \cdot  {\GG}^l_\HH(\xx,\yy,\bsigma)  \\
\ss \cdot  {\GG}^l_\HH(\xx,\yy,\btau) &  \ss \cdot  {\GG}^l_\HH(\xx,\yy,\bsigma)\end{array}\right].
\end{equation}
Each block matrix can then be represented by
\begin{align}
\AA_{0,0} = \tilde{\FF}^0(\xx_i,  \yy_j^{-}, \btau_j^{-}, \bsigma_j^{-}) ~~\mbox{ and }~~ \AA_{0,1} = -\tilde{\FF}^1(\xx_i,  \yy_j^{+}, \btau_j^{+}, \bsigma_j^{+}),
\end{align}
where $i=1,2,\ldots,\cM$ and $j=1,2,\ldots,\cN$,
\begin{align}
\BB_{0,0}&=\FF^0(\xx_i,  \zz_{0,j}, \brho_{0,j},  \bxi_{0,j}) ~~\mbox{ and }~~ \BB_{0,1}=-\FF^1(\xx_i,  \zz_{1,j}, \brho_{1,j},  \bxi_{1,j}),
\end{align}
where $i=1,2,\ldots,\cM$ and $j=1,2,\ldots,\cP$. The right-hand side is the incident fields $\ff = -[\tt_i\cdot \EE^{\inc}, \ss_i\cdot \EE^{\inc}, \mathbf{0}, \mathbf{0} ]^t$, and  the unknown coefficients are arranged as $\cc^- = [\cc^{\btau-}, \cc^{\bsigma-}]^t$, $\cc^+ = [\cc^{\btau+}, \cc^{\bsigma+}]^t$, and $\dd_l = [\dd_l^{\brho}, \dd_l^{\bxi}]^t$ for $l=0,1$. Note that $\AA_{0,0}$ and $\AA_{0,1}$ represent the near interactions of the electric and magnetic fields in the first and second layers, respectively. $\BB_{0,0}$ and $\BB_{0,1}$ represent the far interactions from distant copies of the unit cell in the first and second layers, respectively. 

\paragraph{Quasi-periodicity condition:} The electric and magnetic fields in both the first ($l=0$) and second ($l=1$) layers must satisfy the quasi-periodic condition. In other words, 
\begin{align}
\alpha_x^{-1}\EE_l(\xx_{R_l})-\EE_l(\xx_{L_l}) =0, &~~~~~ \alpha_y^{-1}\EE_l(\xx_{F_l})-\EE_l(\xx_{B_l})=0,\\
\frac{\partial}{\partial \nn}\left(\alpha_x^{-1}\EE_l(\xx_{R_l})-\EE_l(\xx_{L_l}) \right)=0, &~~~~~\frac{\partial}{\partial \nn}\left( \alpha_y^{-1}\EE_l(\xx_{F_l})-\EE_l(\xx_{B_l})\right)=0,
\end{align}
where $\nn$ is the unit normal vector at a target point on one of the side walls. For brevity, only the quasi-periodicity of the electric field for the left and right walls in the first layer is presented. For a point on the left wall $\xx_{L_0,i}$ and the right wall $\xx_{R_0,i} = \xx_{L_0,i}+\ee_x$ in the first layer, the electric field must satisfy
\begin{equation}
\alpha_x^{-1}\EE_0(\xx_{R_0,i})-\EE_0(\xx_{L_0,i})=0.
\end{equation}
By substituting the electric field representation in the first layer and using translation symmetry, this equation can be expanded as
\begin{align}
&\sum_{j=1}^{\cN} c^{\btau-}_{j} \left(\alpha_x^{-2} \sum_{n=-1}^1 \alpha_y^n {\GG}^0_\EE(\xx_{R_0,i},\yy^{-}_{j}-\ee_x+n\ee_y, \btau^-_j)-\alpha_x^{1} \sum_{n=-1}^1 \alpha_y^{n}{\GG}^0_\EE(\xx_{L_0,i}, \yy^{-}_{j}+\ee_x+n\ee_y,\btau^-_j)\right)\nonumber\\
+&\sum_{j=1}^{\cN} c^{\bsigma-}_{j} \left(\alpha_x^{-2} \sum_{n=-1}^1 \alpha_y^n {\GG}^0_\EE(\xx_{R_0,i},\yy^{-}_{j}-\ee_x+n\ee_y, \bsigma^-_j)-\alpha_x^{1} \sum_{n=-1}^1 \alpha_y^{n}{\GG}^0_\EE(\xx_{L_0,i}, \yy^{-}_{j}+\ee_x+n\ee_y,\bsigma^-_j)\right)\nonumber\\
&~~~~~~~~~~+\sum_{j=1}^{\cP}d^{\brho}_{0,j}\left(\alpha_x^{-1}\GG^0_\EE(\xx_{R_0,i},\zz_{0,j},\brho_{0,j})-\GG^0_\EE(\xx_{L_0,i},\zz_{0,j},\brho_{0,j})\right)\nonumber\\
&~~~~~~~~~~+\sum_{j=1}^{\cP}d^{\bxi}_{0,j}\left(\alpha_x^{-1}\GG^0_\EE(\xx_{R_0,i},\zz_{0,j},\bxi_{0,j})-\GG^0_\EE(\xx_{L_0,i},\zz_{0,j},\bxi_{0,j})\right) = 0.
\end{align}
The same procedure can be applied to all other quasi-periodicity conditions and each can be represented in matrix form. Again for notational simplicity, we define 
\small\begin{align}
\tilde{\FF}_Q^l(\xx_w, \yy, \btau):=
 \left[\begin{array}{l}
 \displaystyle\sum_{n=-1}^1 \alpha_y^n\left(\alpha_x^{-2} {\GG}^l_\EE(\xx_{R},\yy-\ee_x+n\ee_y, \btau)-\alpha_x {\GG}^l_\EE(\xx_{L}, \yy+\ee_x+n\ee_y,\btau)\right) \\
 \displaystyle\sum_{n=-1}^1 \alpha_y^n\frac{\partial}{\partial \nn} \left(\alpha_x^{-2} {\GG}^l_\EE(\xx_{R},\yy-\ee_x+n\ee_y, \btau)-\alpha_x {\GG}^l_\EE(\xx_{L}, \yy+\ee_x+n\ee_y,\btau)\right) \\
 \displaystyle\sum_{m=-1}^1 \alpha_x^m \left(\alpha_y^{-2}  {\GG}^l_\EE(\xx_{F},\yy+m\ee_x-\ee_y, \btau)-\alpha_y{\GG}^l_\EE(\xx_{B}, \yy+m\ee_x+\ee_y,\btau)\right) \\
 \displaystyle\sum_{m=-1}^1 \alpha_x^m \frac{\partial}{\partial \nn}\left(\alpha_y^{-2}  {\GG}^l_\EE(\xx_{F},\yy+m\ee_x-\ee_y, \btau)-\alpha_y{\GG}^l_\EE(\xx_{B}, \yy+m\ee_x+\ee_y,\btau)\right) 
 \end{array}\right],
\end{align}\normalsize
and
\begin{align}
&\FF_Q^l(\xx_w, \zz, \brho):= \left[\begin{array}{l}
\alpha_x^{-1}\GG^l_\EE(\xx_{R},\zz,\brho)-\GG^l_\EE(\xx_{L},\zz,\brho) \\
\frac{\partial}{\partial \nn}\left(\alpha_x^{-1}\GG^l_\EE(\xx_{R},\zz,\brho)-\GG^l_\EE(\xx_{L},\zz,\brho)\right) \\
\alpha_y^{-1}\GG^l_\EE(\xx_{F},\zz,\brho)-\GG^l_\EE(\xx_{B},\zz,\brho) \\
\frac{\partial}{\partial \nn}\left(\alpha_y^{-1}\GG^l_\EE(\xx_{F},\zz,\brho)-\GG^l_\EE(\xx_{B},\zz,\brho)\right) \\
\end{array}\right],
\end{align}
where $\xx_{w} = \left( \xx_{R},\xx_{L}, \xx_{F}, \xx_{B}\right)$ denotes $4\cW$ collocation points across all four side walls---$\cW$ points on each wall. For $\{\xx_{w_l,i} = \left( \xx_{R_l,i},\xx_{L_l,i}, \xx_{F_l,i}, \xx_{B_l,i}\right) \}_{i=1}^{\cW}$, with $l=0,1$, the quasi-periodic conditions reduce to 
\begin{equation}
\left[\begin{array}{cc}  \PP_0 &  \bzero \\ \bzero & \PP_1 \end{array}\right]\left[\begin{array}{c}  \cc^-\\ \cc^+ \end{array}\right]+\left[\begin{array}{cc}  \QQ_0 &  \bzero \\ \bzero & \QQ_1 \end{array}\right]\left[\begin{array}{c}  \dd_0\\ \dd_1 \end{array}\right] = \left[\begin{array}{c}  \bzero\\ \bzero\end{array}\right],
\end{equation}
where each block matrix is defined by
\small\begin{equation}
\PP_{0} =\left[\tilde{\FF}_Q^0(\xx_{w_0,i}, \yy^-_j, \btau^-_j)~~ \tilde{\FF}_Q^0(\xx_{w_0,i}, \yy^-_j, \bsigma^-_j)\right] ~~\textrm{and }~~ \PP_{1} =\left[\tilde{\FF}_Q^1(\xx_{w_1,i}, \yy^+_j, \btau^+_j)~~ \tilde{\FF}_Q^1(\xx_{w_1,i}, \yy^+_j, \bsigma^+_j)\right]
\end{equation}\normalsize
for $i=1,2,\ldots,\cW$ and $j=1,2,\ldots,\cN$, and
\begin{align}
\QQ_{l} &=\left[ \FF_Q^l(\xx_{w_l,i}, \zz_{l,j}, \brho_{l,j})~~ \FF_Q^l(\xx_{w_l,i},  \zz_{l,j}, \bxi_{l,j})\right], ~~\textrm{with } l=0,1,
\end{align}
for $i=1,2,\ldots,\cW$ and $j=1,2,\ldots,\cP$.

\paragraph{Radiation condition:} The electric field $\EE_0$ at $\xx_u = (x,y,z_u)$ and $\EE_1$ at $\xx_d = (x,y,z_d)$, and their normal derivatives, must be equal to the Rayleigh-Bloch expansions in Eqs. (\ref{eq:radation1}) and (\ref{eq:radation2}):
\small\begin{align}
\EE_0(\xx_u)=\sum_{m,n \in \bbZ} \aa_u^{mn}\exp(\iu(\kappa_x^m x + \kappa_y^n y)), ~~~&\frac{\partial}{\partial \nn}\EE_0(\xx_u)=\sum_{m,n \in \bbZ} \iu k_{0,z}^{mn}\aa_u^{mn}\exp(\iu(\kappa^m_x x + \kappa_y^n y)), \\
\EE_1(\xx_d)=\sum_{m,n \in \bbZ} \aa_d^{mn}\exp(\iu(\kappa_x^m x + \kappa_y^n y)),~~ ~&\frac{\partial}{\partial \nn}\EE_1(\xx_d)=\sum_{m,n \in \bbZ} \iu k_{1,z}^{mn}\aa_d^{mn}\exp(\iu(\kappa^m_x x + \kappa_y^n y)).
\end{align}\normalsize
Since $\HH$ is coupled with $\EE$ through Maxwell's equations, it suffices to impose the radiation condition only on $\EE$. Expanding the first equation, we find that
\begin{align}
&\sum_{j=1}^{\cN} c^{\btau-}_{j} \tilde{\GG}^0_\EE(\xx_u,\yy^{-}_{j},\btau^-_j) + \sum_{j=1}^{\cN} c^{\bsigma-}_{j} \tilde{\GG}^0_\EE(\xx_u,\yy^{-}_{j},\bsigma^-_j)+\sum_{j=1}^{\cP}d^{\brho}_{0,j}\GG^0_\EE(\xx_u,\zz_{0,j},\brho_{0,j})\nonumber\\
&~~~~~~~~~~~~~~+\sum_{j=1}^{\cP}d^{\bxi}_{0,j}\GG^0_\EE(\xx_u,\zz_{0,j},\bxi_{0,j})-\sum_{m,n \in \bbZ} \aa_u^{mn}\exp(\iu(\kappa_x^m x + \kappa_y^n y)) = \mathbf{0}.
\end{align}
For $\{\xx_{u,i} = (x_i, y_i, z_{u,i})\}_{i=1}^\cW \in \xx_u$ and $\{\xx_{d,i} = (x_i, y_i, z_{d,i})\}_{i=1}^\cW \in \xx_d$, expanding the remaining equations reveals that all four equations can be represented as the matrix equation
\begin{equation}
\left[\begin{array}{cc}\ZZ_0 & 0 \\0 & \ZZ_1\end{array}\right]\left[\begin{array}{c}\cc^- \\ \cc^+\end{array}\right]+\left[\begin{array}{cc}\VV_0 & 0 \\0 & \VV_1\end{array}\right]\left[\begin{array}{c}\dd_0 \\ \dd_1\end{array}\right]+\left[\begin{array}{cc}\WW_0 & 0 \\0 & \WW_1\end{array}\right]\left[\begin{array}{c}\aa_u \\ \aa_d\end{array}\right] = \mathbf{0},
\end{equation}
where the Rayleigh-Bloch expansions have been truncated at $m,n=\pm \cR$, such that $\aa_u = [\aa_u^{mn}]$ and $\aa_d = [\aa_d^{mn}]$ for $m=-\cR,\ldots,\cR$ and $n= -\cR,\ldots,\cR$,
\begin{align}
\ZZ_0 &= \left[\begin{array}{cc}\tilde{\GG}^0_\EE(\xx_{u,i},\yy^{-}_{j},\btau^-_j)  & \tilde{\GG}^0_\EE(\xx_{u,i},\yy^{-}_{j},\bsigma^-_j) \\ \frac{\partial}{\partial \nn}\tilde{\GG}^0_\EE(\xx_{u,i},\yy^{-}_{j},\btau^-_j)  & \frac{\partial}{\partial \nn}\tilde{\GG}^0_\EE(\xx_{u,i},\yy^{-}_{j},\bsigma^-_j) \end{array}\right] ~~\textrm{ and }\\
\ZZ_1 &= \left[\begin{array}{cc}\tilde{\GG}^1_\EE(\xx_{d,i},\yy^{+}_{j},\btau^+_j)  & \tilde{\GG}^1_\EE(\xx_{d,i},\yy^{+}_{j},\bsigma^+_j) \\ \frac{\partial}{\partial \nn}\tilde{\GG}^1_\EE(\xx_{d,i},\yy^{+}_{j},\btau^+_j)  & \frac{\partial}{\partial \nn}\tilde{\GG}^1_\EE(\xx_{d,i},\yy^{+}_{j},\bsigma^+_j) \end{array}\right]
\end{align}
for $i=1,2,\ldots, \cW$ and $j=1,2,\ldots, \cN$,
\begin{align}
\VV_0 &=\left[\begin{array}{cc} \GG^0_\EE(\xx_{u,i},\zz_{0,j},\brho_{0,j}) & \GG^0_\EE(\xx_{u,i},\zz_{0,j},\bxi_{0,j}) \\
\frac{\partial}{\partial \nn}\GG^0_\EE(\xx_{u,i},\zz_{0,j},\brho_{0,j}) & \frac{\partial}{\partial \nn}\GG^0_\EE(\xx_{u,i},\zz_{0,j},\bxi_{0,j})\end{array}\right] ~~\textrm{ and }\\
\VV_1 &=\left[\begin{array}{cc} \GG^1_\EE(\xx_{d,i},\zz_{1,j},\brho_{1,j}) & \GG^1_\EE(\xx_{d,i},\zz_{1,j},\bxi_{1,j}) \\
\frac{\partial}{\partial \nn}\GG^1_\EE(\xx_{d,i},\zz_{1,j},\brho_{1,j}) & \frac{\partial}{\partial \nn}\GG^1_\EE(\xx_{d,i},\zz_{1,j},\bxi_{1,j})\end{array}\right]
\end{align}
for $i=1,2,\ldots, \cW$ and $j=1,2,\ldots, \cP$, and 
\begin{align}
\WW_l = -\left[\begin{array}{ccc}
\texttt{diag}(\exp(\iu(\kappa_x^m x_i + \kappa_y^n y_i)))\\
\texttt{diag}(\iu k_{l,z}^{mn}\exp(\iu(\kappa_x^m x_i + \kappa_y^n y_i)))
\end{array}\right] ~~\textrm{with } l=0,1,
\end{align}
for $i=1,2,\ldots,\cW$, $m=-\cR, \ldots, \cR$, and $n= -\cR, \ldots, \cR$, where $\texttt{diag}(W)$ denotes a $3\times3$ block-diagonal matrix consisting of three copies of $W$ in the diagonal.

\paragraph{Summary:} All three conditions can be combined into the single matrix equation
\begin{equation}
\left[\begin{array}{cccccc} 
\AA_{0,0} & \AA_{0,1}  & \BB_{0,0}  &  \BB_{0,1} & \bzero  & \bzero  \\
\PP_0  &  \bzero &  \QQ_0 &  \bzero & \bzero  & \bzero  \\
\bzero  & \PP_1  & \bzero  & \QQ_1  &  \bzero &  \bzero \\ 
 \ZZ_0 & \bzero  &  \VV_0 &  \bzero &  \WW_0 &  \bzero \\
 \bzero  & \ZZ_1  & \bzero  &  \VV_1 & \bzero  & \WW_1 \end{array}\right]
 \left[\begin{array}{c}  \cc^- \\ \cc^+ \\  \dd_0 \\  \dd_1 \\ \aa_u \\ \aa_d \end{array}\right] = \left[\begin{array}{c} \ff \\ \bzero  \\  \bzero \\  \bzero \\ \bzero \end{array}\right].
\end{equation}
The first row enforces continuity of the tangential components of the electric and magnetic fields; the second and third rows impose the quasi-periodic conditions on the left, right, front, and back walls of the unit cell in the first and second layers; and the last two rows enforce the outgoing radiation condition at both the top and bottom layers. As noted in our Helmholtz equation work \cite{cho2019spectrally}, this is an overdetermined system.  A backward-stable least-squares solver in MATLAB ({\texttt{mldivide}}) can be applied to obtain an accurate solution, or one may use a Schur complement to eliminate $\dd$ and $\aa$ to solve the system.

\begin{figure}[htbp] 
   \centering
   \includegraphics[width=5.2in]{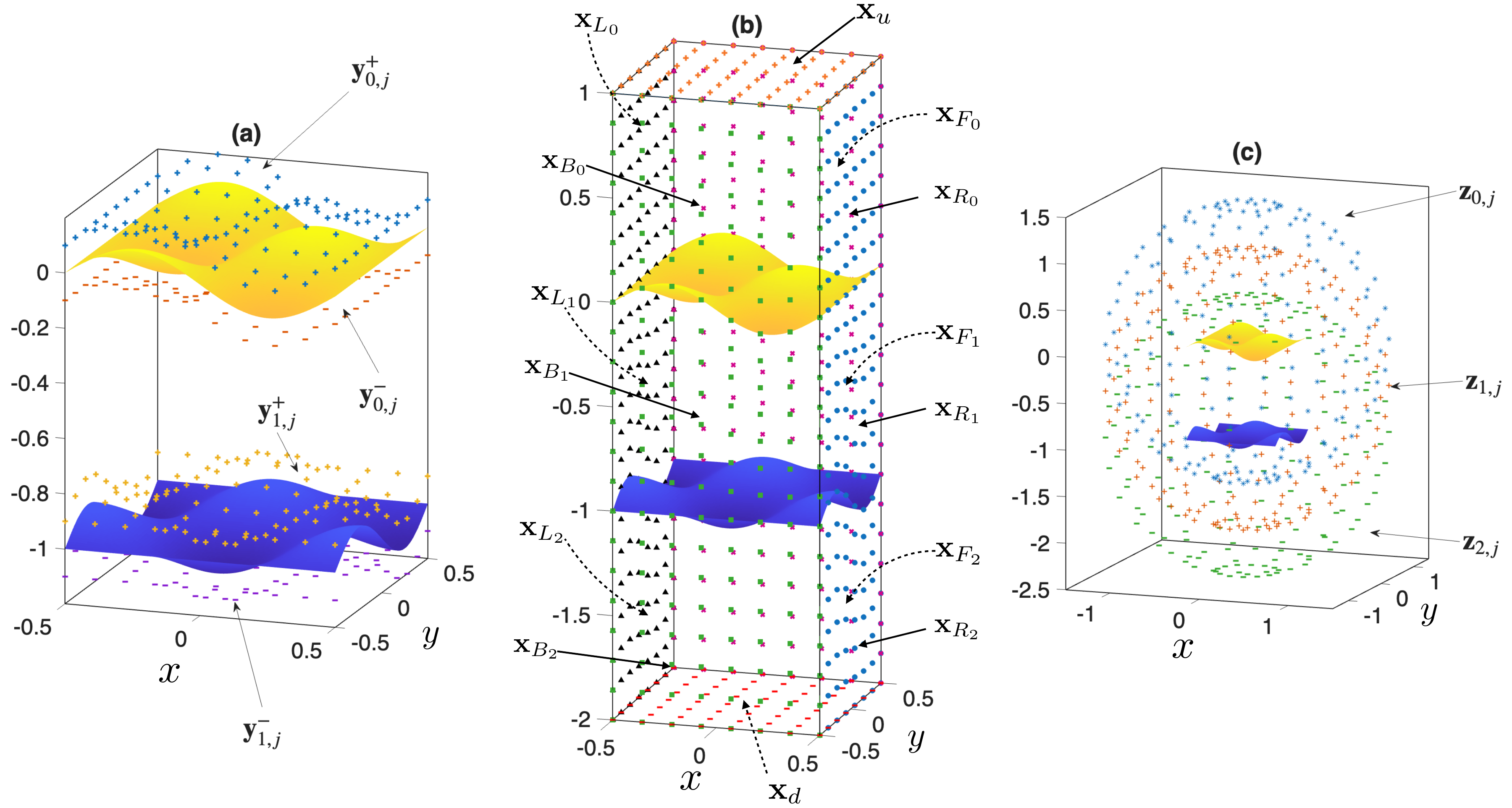} 
   \caption{(a) Two interfaces and MFS source points above and below each interface. (b) The unit cell and surrounding walls for each layer, and (c) proxy source points for each layer.}
   \label{fig:multi_domain}
\end{figure}

\section{MFS for multiple bi-periodic interfaces} \label{sec:multi}
In this section, we generalize the formulation for a single interface to multiple interfaces. The notation extends naturally as follows.
Let $\{\Gamma_l\}_{l=0}^{\cL}$ be bi-periodic interfaces defined by $f_l(x,y)$ with periodicity $d_x$ and $d_y$ along $x$- and $y$-directions. There are $\{\Omega_l\}_{l=0}^{\cL+1}$ domains or layers, where $\Omega_l$ is the region between $f_{l}$ and $f_{l-1}$ surrounded by side walls $\xx_{L_{l}}$, $\xx_{R_{l}}$, $\xx_{B_{l}}$, and $\xx_{F_{l}}$. The first and last layers $\Omega_0$ and $\Omega_{\cL+1}$ have artificial walls at $\xx_u$ and $\xx_d$ for imposing the radiation conditions, respectively. For the periodization scheme, each layer has proxy source points ${\zz_l}$ with two tangent vectors $\brho_l$ and $\bxi_l$ centered on $\Omega_l$ surrounding the unit cell and its immediate eight neighboring cells. Let $\yy_{l}^\pm$ be the coordinate of the MFS source points above (+) and below (-) the interface $f_l$, and $\btau_l^{\pm}$ and $\bsigma_l^{\pm}$ be two tangent vectors at $\yy_{l}^{\pm}$. Figure \ref{fig:multi_domain} shows a notational example for two interfaces with three layers.

For $\xx \in \Omega_0$ (first layer), the MFS with a periodization representation of the electric field is
\begin{align}
    \EE_0(\xx) &\approx  \sum_{j=1}^{\cN} c^{\btau-}_{0,j} \tilde{\GG}^0_\EE(\xx,\yy^{-}_{0,j},\btau^-_{0,j}) + \sum_{j=1}^{\cN} c^{\bsigma-}_{0,j} \tilde{\GG}^0_\EE(\xx,\yy^{-}_{0,j},\bsigma^-_{0,j})\nonumber\\
    &+\sum_{j=1}^{\cP}d^{\brho}_{0,j}\GG^0_\EE(\xx,\zz_{0,j},\brho_{0,j})+\sum_{j=1}^{\cP}d^{\bxi}_{0,j}\GG^0_\EE(\xx,\zz_{0,j},\bxi_{0,j}).\label{eq:mlp1} 
  \end{align}
For $\xx\in \{\Omega_l\}_{l=1}^{\cL}$ (interior layers), the interfaces above and below contribute to the electric fields such that
    \begin{align}
    \EE_l(\xx) &\approx \sum_{j=1}^{\cN} c^{\btau+}_{l-1,j} \tilde{\GG}^l_\EE(\xx,\yy^{+}_{l-1,j},\btau^+_{l-1,j}) + \sum_{j=1}^{\cN} c^{\bsigma+}_{l-1,j} \tilde{\GG}^l_\EE(\xx,\yy^{+}_{l-1,j},\bsigma^+_{l-1,j})\nonumber\\
 &+   \sum_{j=1}^{\cN} c^{\btau-}_{l,j} \tilde{\GG}^l_\EE(\xx,\yy^{-}_{l,j},\btau^-_{l,j}) + \sum_{j=1}^{\cN} c^{\bsigma-}_{l,j} \tilde{\GG}^l_\EE(\xx,\yy^{-}_{l,j},\bsigma^-_{l,j})\nonumber\\
 &~~~~~~~~~~~~~~~~~~~~~~~~~+\sum_{j=1}^{\cP}d^{\brho}_{l,j}\GG^l_\EE(\xx,\zz_{l,j},\brho_{l,j})+\sum_{j=1}^{\cP}d^{\bxi}_{l,j}\GG^l_\EE(\xx,\zz_{l,j},\bxi_{l,j}).\label{eq:mlp2} 
    \end{align}
For $\xx \in \Omega_{\cL+1}$ (last layer), only the interface MFS source points above $\Gamma_l$ contribute to the electric field such that
    \begin{align}
\EE_{\cL+1}(\xx) &\approx  \sum_{j=1}^{\cN} c^{\btau+}_{{\cL},j} \tilde{\GG}^{\cL+1}_\EE(\xx,\yy^{+}_{\cL, j},\btau^+_{\cL,j}) + \sum_{j=1}^{\cN} c^{\bsigma+}_{{\cL},j} \tilde{\GG}^{\cL+1}_\EE(\xx,\yy^{+}_{\cL, j},\bsigma^+_{\cL, j}) \nonumber\\
    &+\sum_{j=1}^{\cP}d^{\brho}_{{\cL+1},j}\GG^{\cL+1}_\EE(\xx,\zz_{\cL+1,j},\brho_{\cL+1,j})+\sum_{j=1}^{\cP}d^{\bxi}_{{\cL+1},j}\GG^{\cL+1}_\EE(\xx,\zz_{\cL+1,j},\bxi_{\cL+1,j}).\label{eq:mlp3} 
 \end{align}
By symmetry, the formulas for the magnetic fields can be expressed by substituting $\EE$ by $\HH$, and have therefore been omitted. The unknown coefficients $\cc_l^{\tau\pm} = \{c_{l,j}^{\tau\pm} \}_{j=1}^{\cN}$, $\cc_l^{\sigma\pm} = \{c_{l,j}^{\sigma\pm} \}_{j=1}^{\cN}$, $\dd_l^{\rho} =\{d_{l,j}^{\rho }\}_{j=1}^{\cP}$, and $\dd_l^{\xi} =\{d_{l,j}^{\xi }\}_{j=1}^{\cP}$ can be found by enforcing the continuity of the tangential components of the electric and magnetic fields at each interface $\{ \Gamma_l\}_{l=0}^{\cL}$, quasi-periodicity on the side walls of each $\Omega_l$, and Rayleigh-Bloch expansions in the topmost  and bottommost layers.

\paragraph{Transmission condition:} For the incident wave $\EE^{\inc}(\xx) = \EE^{\inc}e^{\iu \kk \cdot \xx}$ in $\Omega_0$, the continuity of the tangential components is
\begin{align}
\tt \cdot \EE_0(\xx) - \tt \cdot \EE_{1}(\xx) &= -\tt\cdot \EE^{\inc}(\xx),\\
\ss \cdot \EE_0(\xx) - \ss \cdot \EE_{1}(\xx) &= -\ss\cdot \EE^{\inc}(\xx),\\
\tt \cdot \HH_0(\xx) - \tt \cdot \HH_{1}(\xx) &= \mathbf{0},\\
\ss \cdot \HH_0(\xx) - \ss \cdot \HH_{1}(\xx) &= \mathbf{0},
\end{align}
for $\xx \in \Gamma_0$, and 
\begin{align}
\tt \cdot \EE_l(\xx) - \tt \cdot \EE_{l+1}(\xx) &= \mathbf{0},\\
\ss \cdot \EE_l(\xx) - \ss \cdot \EE_{l+1}(\xx) &= \mathbf{0},\\
\tt \cdot \HH_l(\xx) - \tt \cdot \HH_{l+1}(\xx) &= \mathbf{0},\\
\ss \cdot \HH_l(\xx) - \ss \cdot \HH_{l+1}(\xx) &= \mathbf{0},
\end{align}
for $\xx\in \{\Gamma_l\}_{l=1}^{\cL-1}$.

Let the unknown coefficients be $\cc_0^{-} = [\cc_{0}^{\btau-}, \cc_{0}^{\bsigma-}]^t$, $\cc_l = [\cc_{l-1}^{\btau+}, \cc_{l-1}^{\bsigma+}, \cc_l^{\btau-}, \cc_l^{\bsigma-}]^t$ for $l=1,2,\ldots,\cL$, $\cc_\cL^{+} = [\cc_{\cL}^{\btau+}, \cc_{\cL}^{\bsigma+}]^t$, and $\dd_l = [\dd_{l}^{\brho}, \dd_{l}^{\bxi}]^t$, for  $l=0,1,\ldots,\cL+1$. All the continuity conditions can then be represented as the matrix equation
 \begin{equation}
 \AA \cc+\BB \dd = \ff,
 \end{equation}
 where
  \small\begin{equation}
 \AA= \left[\begin{array}{cccccc} 
  \AA_{0,0} &\AA_{0,1} &\bzero                &\cdots  &\bzero  \\
   \bzero    &\AA_{1,1}  &\AA_{1,2}           &\cdots  &\bzero \\  
   \vdots    &   \vdots       &   \vdots               &\vdots  &\vdots \\  
   \bzero    &     \bzero       &      \cdots         &\AA_{\cL, \cL}  &\AA_{\cL, \cL+1}   
    \end{array}\right],
 ~~\BB= \left[\begin{array}{cccccc} 
\BB_{0,0} &\BB_{0,1} & \bzero  & \cdots & \bzero\\
\bzero      &\BB_{1,1} & \BB_{1,2}& \cdots & \bzero\\  
\vdots      &   \vdots   &    \vdots    & \vdots & \vdots\\  
\bzero      &   \bzero   &    \cdots   & \BB_{\cL, \cL} & \BB_{\cL, \cL+1}\\  
    \end{array}\right],
   \end{equation}\normalsize
$\ff = [-\tt_i\cdot \EE^{\inc}(\xx_i), -\ss_i\cdot \EE^{\inc}(\xx_i), \bzero,\ldots,\bzero ]^t$ representing the incident field impinging on the interface $\Gamma_0$, and the unknown vectors $\cc$ and $\dd$ are arranged as
$$\cc= \left[\begin{array}{c} 
\cc_0^-, 
\cc_1,
\cc_2,
\cdots,
\cc_\cL,
\cc_\cL^+
\end{array}\right] ^t \mbox{ and }
\dd= \left[\begin{array}{c} 
\dd_{0},
\dd_{1},
\dd_{2},
\cdots,
\dd_{\cL+1}\\
\end{array}\right] ^t.
$$
The diagonal blocks $\AA_{l,l}$ and $\BB_{l,l}$ represent contributions from the layer $\Omega_l$, and the off-diagonal blocks $\AA_{l,l+1}$ and $\BB_{l,l+1}$ represent contributions from the layer $\Omega_{l+1}$. Therefore, each row enforces the continuity condition at the interface $\Gamma_l$. Note that the first and last block matrices of $\AA$ have only one interface contributing to their solutions and are given by
\begin{align}
&\AA_{0,0} =\tilde{\FF}^0(\xx_i,  \yy_{0,j}^-, \btau_{0,j}^{-}, \bsigma_{0,j}^-) ~~\textrm{and}~~ \AA_{\cL,\cL+1} =-\tilde{\FF}^{\cL+1}(\xx_i,  \yy_{\cL,j}^+, \btau_{\cL,j}^{+}, \bsigma_{\cL,j}^+ ).
\end{align}
The remaining diagonal block matrices for $l=1,2,\ldots,\cL$ are given by
\begin{align}
&\AA_{l,l} =\left[\begin{array}{cc}\tilde{\FF}^l(\xx_i,  \yy_{l-1,j}^+, \btau_{l-1,j}^{+}, \bsigma_{l-1,j}^+) & \tilde{\FF}^l(\xx_i,   \yy_{l,j}^-, \btau_{l,j}^{-}, \bsigma_{l,j}^-)\end{array}\right].
\end{align}
The off-diagonal blocks for $l=0,1,\cdots,\cL-1$ are given by
\begin{align}
&\AA_{l,l+1} =-\left[\begin{array}{cc}\tilde{\FF}^{l+1}(\xx_i,  \yy_{l,j}^+, \btau_{l,j}^{+}, \bsigma_{l,j}^+) & \tilde{\FF}^{l+1}(\xx_i,  \yy_{l+1,j}^-, \btau_{l+1,j}^{-}, \bsigma_{l+1,j}^-)\end{array}\right].
\end{align}
Lastly, the contribution from proxy source points with $l=0,1,\ldots,\cL$ are given by
\begin{align}
\BB_{l,l} = \FF^l(\xx_i,  \zz_{l,j}, \brho_{l,j}, \bxi_{l,j}) ~~\mbox{and}~~ \BB_{l,l+1} = -\FF^{l+1}(\xx_i, \zz_{l+1,j}, \brho_{l+1,j}, \bxi_{l+1,j}).
\end{align}

\paragraph{Quasi-periodicity condition:} The electric field in each layer must satisfy the quasi-periodic condition. In other words, 
\begin{align}
\alpha_x^{-1}\EE_l(\xx_{R_l})-\EE_l(\xx_{L_l}) =0, &~~~~~ \alpha_y^{-1}\EE_l(\xx_{F_l})-\EE_l(\xx_{B_l})=0,\\
\frac{\partial}{\partial \nn}\left(\alpha_x^{-1}\EE_l(\xx_{R_l})-\EE_l(\xx_{L_l}) \right)=0, &~~~~~\frac{\partial}{\partial \nn}\left( \alpha_y^{-1}\EE_l(\xx_{F_l})-\EE_l(\xx_{B_l})\right)=0,
\end{align}
where $l=0,1,\ldots, \cL+1$, and $\nn$ is the unit normal vector at a target point on one of the side walls. The quasi-periodicity condition can be represented as the matrix equation
\begin{equation}
\PP \cc+\QQ \dd = \bzero,
\end{equation}
where $\PP$ and $\QQ$ are diagonal block matrices defined by
\begin{equation}
\PP = \left[\begin{array}{ccccc}  
\PP_0 &  \bzero  & \bzero & \cdots & \bzero\\ 
\bzero & \PP_1  & \bzero  & \cdots & \bzero \\
\vdots & \vdots  & \vdots  & \vdots & \vdots \\
\bzero & \bzero  & \bzero  & \cdots & \PP_{\cL+1}
\end{array}\right] ~~\mbox{and}~~
\QQ=\left[\begin{array}{ccccc}  
\QQ_0 &  \bzero  & \bzero & \cdots & \bzero\\ 
\bzero & \QQ_1  & \bzero  & \cdots & \bzero \\
\vdots & \vdots  & \vdots  & \vdots & \vdots \\
\bzero & \bzero  & \bzero  & \cdots & \QQ_{\cL+1}
\end{array}\right].
\end{equation}
Each $\PP_l$ and $\QQ_l$ impose the quasi-periodicity condition for layer $\Omega_l$, and each block is given by
\begin{align}
\PP_{0} &=\left[\tilde{\FF}_Q^0(\xx_{w_0,i}, \yy^-_{0, j}, \btau^-_{0, j})~~ \tilde{\FF}_Q^0(\xx_{w_0,i}, \yy^-_{0, j}, \bsigma^-_{0, j})\right],\\
\PP_{l} &=\left[\tilde{\FF}_Q^l(\xx_{w_l,i}, \yy^-_{l, j}, \btau^-_{l, j})~~ \tilde{\FF}_Q^l(\xx_{w_l,i}, \yy^-_{l, j}, \bsigma^-_{l, j}) \right. \nonumber\\
&~~~~~~~~~~\left.~~ \tilde{\FF}_Q^l(\xx_{w_l,i}, \yy^+_{l-1, j}, \btau^+_{l-1, j})~~ \tilde{\FF}_Q^l(\xx_{w_l,i}, \yy^+_{l-1, j}, \bsigma^+_{l-1, j})\right] ~~\textrm{with } l=1,2,\ldots, \cL,\\
\PP_{\cL+1} &=\left[\tilde{\FF}_Q^{\cL+1}(\xx_{w_{\cL+1},i}, \yy^+_{\cL, j}, \btau^+_{\cL, j})~~ \tilde{\FF}_Q^{\cL+1}(\xx_{w_{\cL+1},i}, \yy^+_{\cL, j}, \bsigma^+_{\cL, j})\right],
\end{align}
for $i=1,2,\ldots, \cW$ and $j=1,2,\ldots, \cN$. Note that $\PP_0$ and $\PP_{\cL+1}$ have only one interface contributing to the electric field in $\Omega_0$ and $\Omega_{\cL+1}$, respectively. The definition of $\QQ$ is the same as the single interface one
\begin{align}
\QQ_{l} &=\left[ \FF_Q^l(\xx_{w_l,i}, \zz_{l,j}, \brho_{l,j})~~ \FF_Q^l(\xx_{w_l,i},  \zz_{l,j}, \bxi_{l,j})\right] ~~ \textrm{with } l=0,1, \ldots, \cL+1,
\end{align}
for $i=1,2,\ldots,\cW$ and $j=1,2,\ldots,\cP$.

\paragraph{Radiation conditions:} Radiation conditions only need to be applied in the topmost and bottommost layers. Therefore, the single-interface formula can be extended by including zero block matrices to separate the first and  last equations and by changing the indices appropriately, such that
\begin{equation}
\ZZ\cc+\VV \dd+\WW\aa = \bzero,
\end{equation}
where
\begin{equation}
\ZZ = \left[\begin{array}{cccc}
\ZZ_0 & \bzero&\cdots &\bzero \\
\bzero &\bzero & \cdots & \ZZ_{\cL+1}\end{array}\right],
~~\VV=\left[\begin{array}{cccc}
\VV_0 & \bzero&\cdots &\bzero \\
\bzero &\bzero & \cdots & \VV_{\cL+1}\end{array}\right],
~~\WW = \left[\begin{array}{cc}
\WW_0 & \bzero \\
\bzero &  \WW_{\cL+1}\end{array}\right],
\end{equation}
\begin{align}
\ZZ_0 &= \left[\begin{array}{cc}\tilde{\GG}^0_\EE(\xx_{u,i},\yy^{-}_{0,j},\btau^-_{0,j})  & \tilde{\GG}^0_\EE(\xx_{u,i},\yy^{-}_{0,j},\bsigma^-_{0,j}) \\ \frac{\partial}{\partial \nn}\tilde{\GG}^0_\EE(\xx_{u,i},\yy^{-}_{0,j},\btau^-_{0,j})  & \frac{\partial}{\partial \nn}\tilde{\GG}^0_\EE(\xx_{u,i},\yy^{-}_{0,j},\bsigma^-_{0,j}) \end{array}\right] ~~\textrm{and} \\
\ZZ_{\cL+1} &= \left[\begin{array}{cc}\tilde{\GG}^{\cL+1}_\EE(\xx_{d,i},\yy^{+}_{\cL,j},\btau^+_{\cL, j})  & \tilde{\GG}^{\cL+1}_\EE(\xx_{d,i},\yy^{+}_{\cL, j},\bsigma^+_{\cL, j}) \\ \frac{\partial}{\partial \nn}\tilde{\GG}^{\cL+1}_\EE(\xx_{d,i},\yy^{+}_{\cL, j},\btau^+_j)  & \frac{\partial}{\partial \nn}\tilde{\GG}^{\cL+1}_\EE(\xx_{d,i},\yy^{+}_{\cL, j},\bsigma^+_{\cL, j}) \end{array}\right]
\end{align}
for $i=1,2,\ldots,\cW$ and $j=1,2,\ldots,\cN$,
\begin{align}
\VV_0 &=\left[\begin{array}{cc} \GG^0_\EE(\xx_{u,i},\zz_{0,j},\brho_{0,j}) & \GG^0_\EE(\xx_{u,i},\zz_{0,j},\bxi_{0,j}) \\
\frac{\partial}{\partial \nn}\GG^0_\EE(\xx_{u,i},\zz_{0,j},\brho_{0,j}) & \frac{\partial}{\partial \nn}\GG^0_\EE(\xx_{u,i},\zz_{0,j},\bxi_{0,j})\end{array}\right] ~~\textrm{and} \\
\VV_{\cL+1} &=\left[\begin{array}{cc} \GG^{\cL+1}_\EE(\xx_{d,i},\zz_{{\cL+1},j},\brho_{{\cL+1},j}) & \GG^{\cL+1}_\EE(\xx_{d,i},\zz_{{\cL+1},j},\bxi_{{\cL+1},j}) \\
\frac{\partial}{\partial \nn}\GG^{\cL+1}_\EE(\xx_{d,i},\zz_{{\cL+1},j},\brho_{{\cL+1},j}) & \frac{\partial}{\partial \nn}\GG^{\cL+1}_\EE(\xx_{d,i},\zz_{{\cL+1},j},\bxi_{{\cL+1},j})\end{array}\right]
\end{align}
for $i=1,2,\ldots,\cW$ and $j=1,2,\ldots,\cP$, and
\begin{align}
\WW_l = -\left[\begin{array}{ccc}
\texttt{diag}(\exp(\iu(\kappa_x^m x_i + \kappa_y^n y_i)))\\
\texttt{diag}(\iu k_{l,z}^{mn}\exp(\iu(\kappa_x^m x_i + \kappa_y^n y_i)))
\end{array}\right] ~~\textrm{with } l=0,\cL+1,
\end{align}
for $i=1,2,\ldots,\cW$, $m=-\cR, \ldots, \cR$, and $n= -\cR, \ldots, \cR$, where $\texttt{diag}(W)$ denotes a $3\times3$ block-diagonal matrix consisting of three copies of $W$ in the diagonal.

\paragraph{Summary:} All three conditions can be combined into a single matrix equation
\begin{align}
\left[\begin{array}{ccc}
\AA & \BB & \bzero \\
\PP & \QQ &\bzero \\
\ZZ & \VV & \WW\end{array}\right]
\left[\begin{array}{c}\cc \\ \dd \\ \aa\end{array}\right]
=\left[\begin{array}{c}\ff \\ \bzero \\ \bzero\end{array}\right],
\end{align}
where $\AA$ and $\BB$ are block bi-diagonal matrices imposing the continuity conditions at each interface, $\PP$ and $\QQ$ are block diagonal matrices for the quasi-periodic conditions, and $\ZZ$, $\VV$, and $\WW$ enforce the Rayleigh-Bloch radiation conditions. The same solution procedure from the Helmholtz equation can be applied. Schur complements eliminate additional unknowns from the periodization scheme, and the resulting matrix has the same size and structure as $\AA$. All of the MFS coefficients $\cc$ can then be found using the backward-stable least-squares solver in MATLAB ({\texttt{mldivide}}). 
 
\begin{figure}[h] 
   \centering
   \includegraphics[width=4.5in]{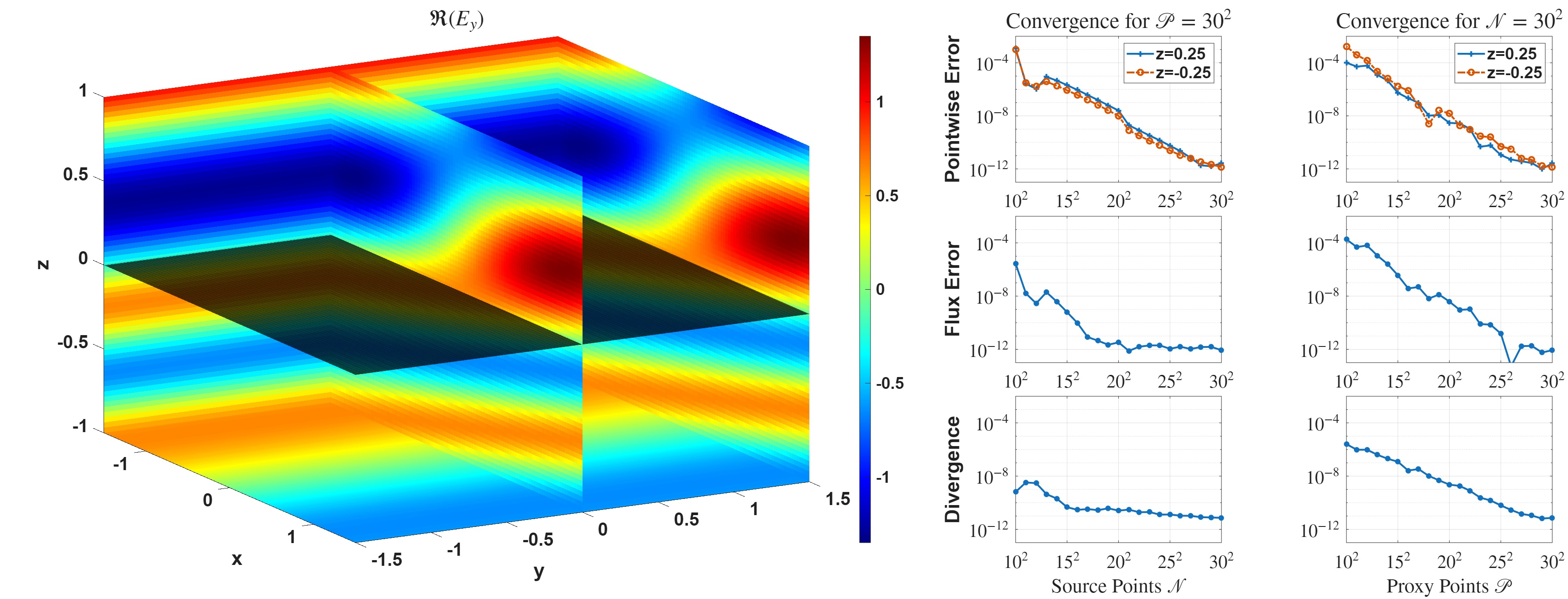} 
   \caption{One flat interface. \textbf{Left:} Real part of $E_y$. \textbf{Center:} Pointwise convergence, flux error, and divergence at $(0,0,0.25)$ ($+$ marker) and $(0,0, -0.25)$ (o marker) as a function of MFS source points against the exact solution. \textbf{Right:} Pointwise convergence, flux error, and divergence at $(0,0,0.25)$ ($+$ marker) and $(0,0, -0.25)$ (o marker) as a function of proxy source points against the exact solution.}
   \label{fig:oneflat}
\end{figure}
 
\section{Numerical Results}\label{sec:experiments}
In this section, we present numerical experiments  conducted using MATLAB R2024a on a dual 2.20 GHz Intel Xeon Gold 5220R CPU (48 cores total) workstation equipped with 768 GB of memory. MATLAB code used to produce graphics and numerical results presented in this paper can be found at \href{https://github.com/jaredweed1990/mfs_multilayer_media}{https://github.com/jaredweed1990/mfs\_multilayer\_media}. Pointwise convergence, divergence-free conditions, and flux error are used to assess performance. The number of target points on each interface is set to $\lceil(1.1)^2\cN\rceil$ for all numerical examples. 

The energy is defined by integrating the normal component of the Poynting vector over a unit surface as
\begin{equation}
\mathcal{E} = \int_S \frac{1}{2} \left(\EE \times \HH^*\right)\cdot \nn~ dS.
\end{equation}
From Ref. \cite{liu2016numerical}, the incident energy ($\mathcal{E}^i$) is derived from the incident field, and the reflected ($\mathcal{E}^r$) and transmitted energy ($\mathcal{E}^t$) are numerically computed using Rayleigh-Bloch expansions with the numerically computed coefficients $\aa_u$ and $\aa_d$:
\begin{align}
\mathcal{E}^i  &= -\frac{1}{2}\frac{d_x d_y}{\omega \mu_0}\cdot \Re\left(k_z {(E_x^{\inc})}^2-k_xE_x^{\inc} E_z^{{\inc},*}-k_yE_y^{{\inc}}E_z^{{\inc},*}+k_z (E_y^{inc})^2 \right),\\
\mathcal{E}^+ &= \frac{d_xd_y}{2\omega \mu_0} \cdot \Re\left(\sum_{m,n=-\cR}^\cR a_{u,x}^{mn}(\kappa_{0,z}^{mn} a_{u,x}^{mn,*}-\kappa_x^{m}a_{u,z}^{mn,*})-a_{u,y}^{mn}(\kappa_{y}^{n} a_{u,z}^{mn,*}-\kappa_{0,z}^{mn}a_{u,y}^{mn,*}) \right),\\
\mathcal{E}^-  &=\frac{d_xd_y}{2\omega \mu_l} \cdot \Re\left(\sum_{m,n=-\cR}^\cR a_{d,x}^{mn}(\kappa_{l,z}^{mn} a_{d,x}^{mn,*}+\kappa_x^{m}a_{d,z}^{mn,*})+a_{d,y}^{mn}(\kappa_{y}^{n} a_{d,z}^{mn,*}+\kappa_{l,z}^{mn}a_{d,y}^{mn,*}) \right).
\end{align}
The flux error is then defined as  $|\mathcal{E}^r+\mathcal{E}^t-\mathcal{E}^i|$. In addition, the divergence condition is computed via a random sampling of interior points, away from cell boundaries and MFS source points, using a standard six-point stencil method.

\paragraph{Error convergence for a flat interface:} For the first experiment, we evaluate our numerical solution at $(0,0,0.25)$ and $(0,0,-0.25)$ and compare it with the exact Fresnel solution \cite{chew1999waves} ($-0.203379977935232 -0.285039015281348\iu$ and $-0.256161363870884 +0.597223512314425\iu$, respectively) for a flat surface $f(x,y) = 0$, where $\epsilon_0 = 1$, $\mu_0=1$, $\epsilon_1 = 4$, $\mu_1=1$, and $\omega = 4$ (see Fig. \ref{fig:oneflat}). The transverse electric (TE) incident field $\EE^{\inc}(\xx) = (0,1,0)e^{\iu \kk \cdot \xx}$ with a wavevector of $\kk$ pointing in the direction of azimuthal angle $\phi = 9\pi/10$ and polar angle $\theta = 0$ and the real part of the electric field $E_y$ for the entire cell is shown. The pointwise convergence, flux error, and divergence are plotted as a function of MFS source points and proxy source points. The number of points $\cW$ on the walls surrounding the unit cell as well as on the top and bottom artificial radiation surfaces remain constant at $\cW=30^2$, with the Rayleigh-Bloch expansions truncated at $\cR=10$. For the convergence with respect to MFS source points, $\cP=30^2$ proxy source points are used. For the convergence with respect to proxy source points, $\cN=30^2$ MFS source points are used. The MFS source points are offset 0.15 units above and below the surface along the normal direction.

\begin{figure}[h] 
   \centering
   \includegraphics[width=4.5in]{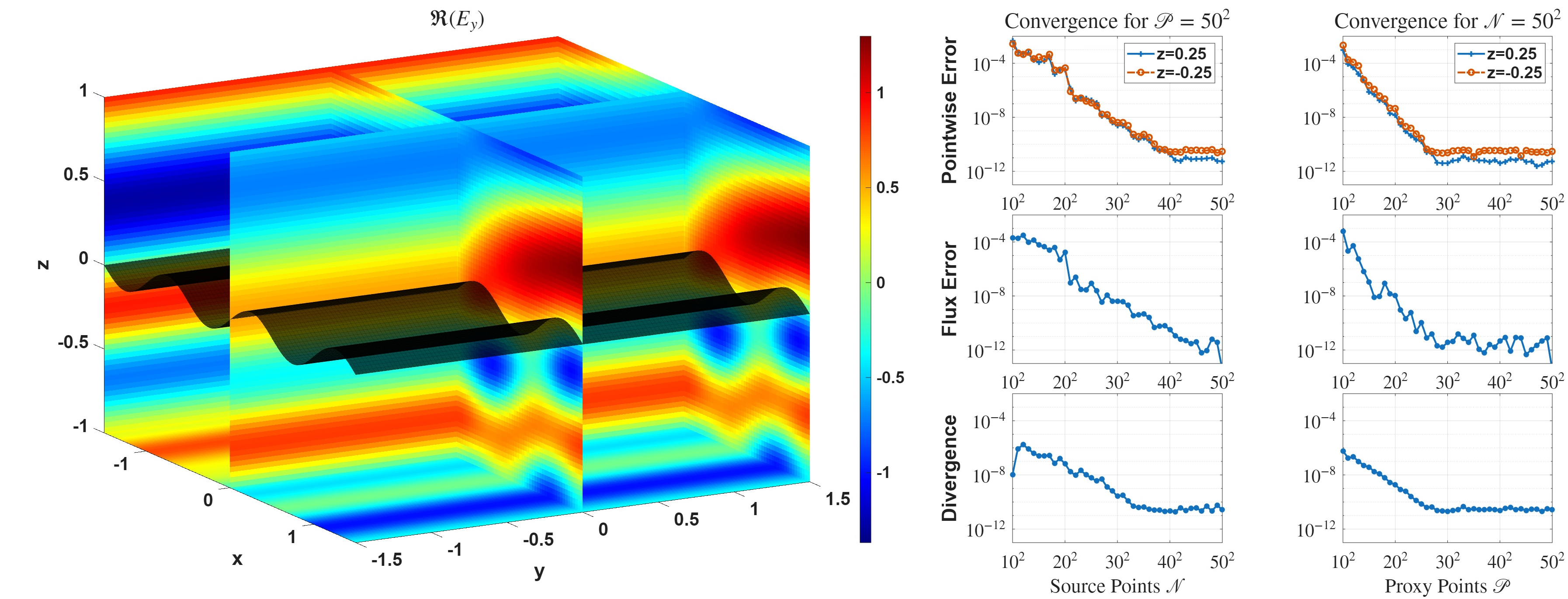} 
   \caption{One interface with $f(x,y) = 0.1\sin(2\pi x)$. \textbf{Left:} Real part of $E_y$. \textbf{Center:} Pointwise convergence, flux error, and divergence at $(0,0,0.25)$ ($+$ marker) and $(0,0, -0.25)$ (o marker) as a function of MFS source points against the reference solution. \textbf{Right:} Pointwise convergence, flux error, and divergence at $(0,0,0.25)$ ($+$ marker) and $(0,0, -0.25)$ (o marker)as a function of proxy source points against the reference solution.}
   \label{fig:onesine}
\end{figure}

\paragraph{Error convergence for a 1D-periodic interface:} For the second experiment, we evaluate our numerical solution at $(0,0,0.25)$ and $(0,0,-0.25)$ and compare it with a reference solution obtained via the Helmholtz equation solver \cite{cho2019spectrally} ($-0.300249349648672 - 0.210359576389989\iu$ and $-0.204763042628003 + 0.273746531105035\iu$, respectively) for a surface $f(x,y) = 0.1\sin(2\pi x)$, where $\epsilon_0=1$, $\mu_0=1$, $\epsilon_1=4$, $\mu_1=1$, and $\omega = 4$ (see Fig. \ref{fig:onesine}). We use the same TE incident field and wavevector $\kk$ here as in the flat interface case. With this setup, we note that Maxwell's equations reduce to the Helmholtz equation, and reference solutions are obtained from the authors' previous work (Ref.\cite{cho2019spectrally}). Figure \ref{fig:onesine} displays the real part of the electric field component $E_y$, the pointwise convergence, flux error, and divergence as functions of both the MFS and proxy source points counts. $\cP=50^2$ proxy points are used for MFS source point convergence, $\cN=50^2$ source points are used for proxy point convergence, and $\cW=30^2$ and $\cR=10$ remain constant. The MFS source points are offset 0.15 units above and below the surface along the surface normal direction.

\begin{figure}[h] 
   \centering
   \includegraphics[width=4.5in]{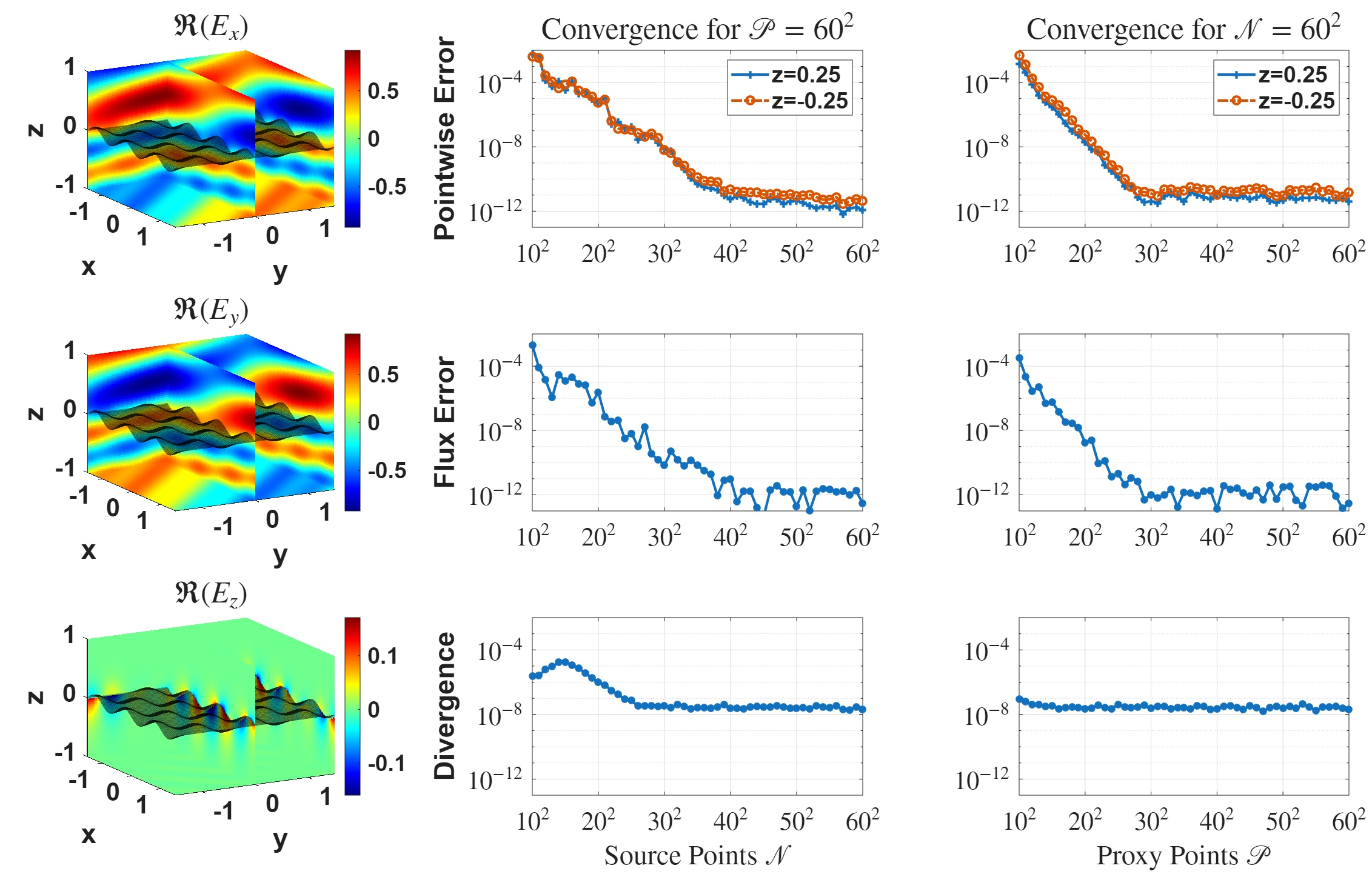} 
   \caption{One interface with $f(x,y) = 0.1\sin(2\pi x)\cos(2\pi y)$. \textbf{Left:} Real part of $E_x$, $E_y$, and $E_z$. \textbf{Center:} Pointwise convergence  at $(0,0,\pm 0.25)$, flux error, and divergence as a function of MFS source points. \textbf{Right:} Pointwise convergence at $(0,0,\pm0.25)$ , flux error, and divergence as a function of proxy source points.}
   \label{fig:onesinecos}
\end{figure}

\paragraph{Error convergence for a bi-periodic interface:} For the third numerical experiment, we define the surface as $f(x,y) = 0.1\sin(2\pi x)\cos(2\pi y)$, where $\epsilon_0=1$, $\mu_0=1$, $\epsilon_1=4$, $\mu_1=1$, and $\omega = 4$. The solution is evaluated at $(0,0,0.25)$ and $(0,0,-0.25)$. The TE incident field $\EE^{\inc}(\xx) = (-\frac{\sqrt{2}}{2},\frac{\sqrt{2}}{2}, 0)e^{\iu \kk \cdot \xx}$ with a wavevector of $\kk$ pointing in the direction of azimuthal angle $\phi = 9\pi/10$ and polar angle $\theta = \pi/4$ and the real part of each component to the electric field for the entire cell is shown, along with pointwise convergence in $E_y$, flux error, and divergence are presented in Figure \ref{fig:onesinecos}. $\cP=60^2$ proxy points are used for MFS source point convergence, $\cN=60^2$ source points are used for proxy point convergence, and $\cW=30^2$ and $\cR=10$ are fixed. The MFS source points are offset 0.15 units above and below the surface along the surface normal.

\begin{figure}[h] 
   \centering
   \includegraphics[width=3in]{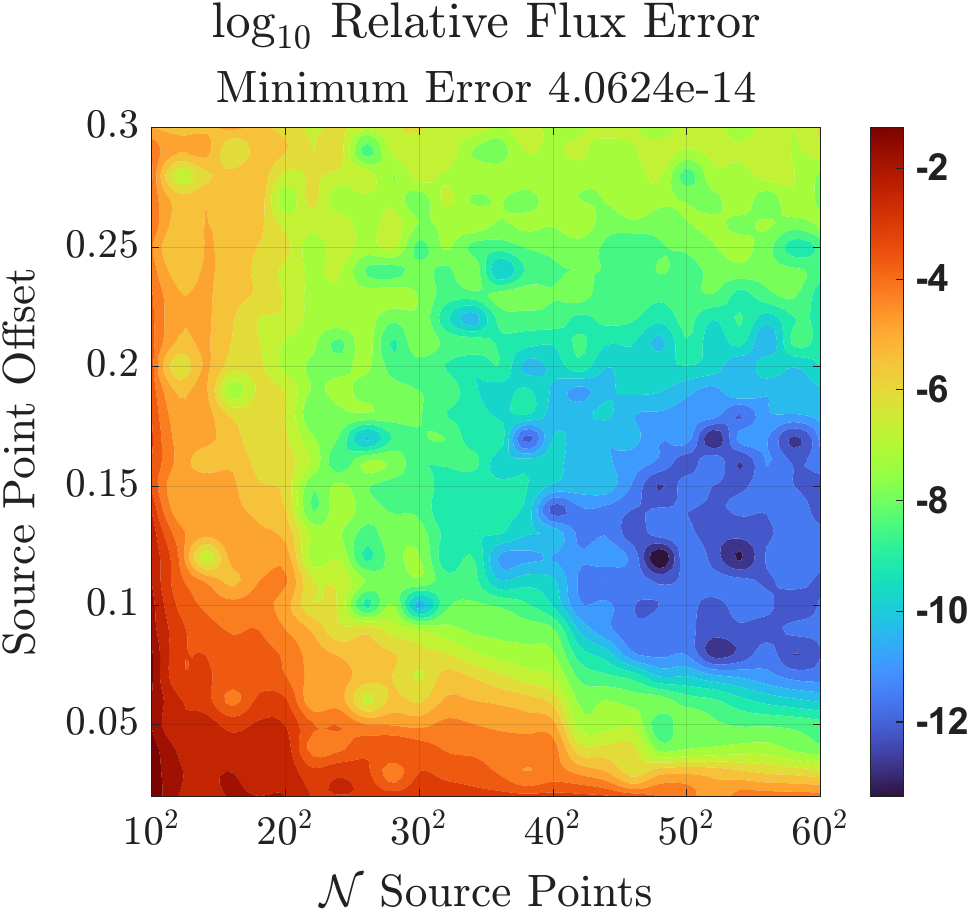} 
   \caption{Flux Error of MFS solutions for a single interface defined as $f(x,y) = 0.1\sin(2\pi x)\cos(2\pi y)$ where the source point offset and number of source points are varied. The minimum error is $4.06\times 10^{-14}$ when $\cN=48^2$ with the offset at $0.12$ units.}
   \label{fig:mfsdistvary}
\end{figure}

\paragraph{MFS source point offset optimality for a bi-periodic interface:} For the fourth numerical experiment, we use $f(x,y)=0.1\sin(2\pi x)\cos(2\pi y)$ to define the surface, where $\epsilon_0=1$, $\mu_0=1$, $\epsilon_1=4$, $\mu_1=1$, and $\omega = 4$, with the TE incident field $\EE^{\inc}(\xx) = (-\frac{\sqrt{2}}{2},\frac{\sqrt{2}}{2}, 0)e^{\iu \kk \cdot \xx}$ and a wavevector of $\kk$ pointing in the direction of azimuthal angle $\phi = 9\pi/10$ and polar angle $\theta = \pi/4$. The flux error as a function of MFS source points offsets and total source points $\cN$ used is presented in Figure \ref{fig:mfsdistvary}. For this experiment, $\cR=10$, $\cW=30^2$ and $\cP=40^2$ proxy points are fixed.

\begin{figure}[h] 
   \centering
   \includegraphics[width=4.5in]{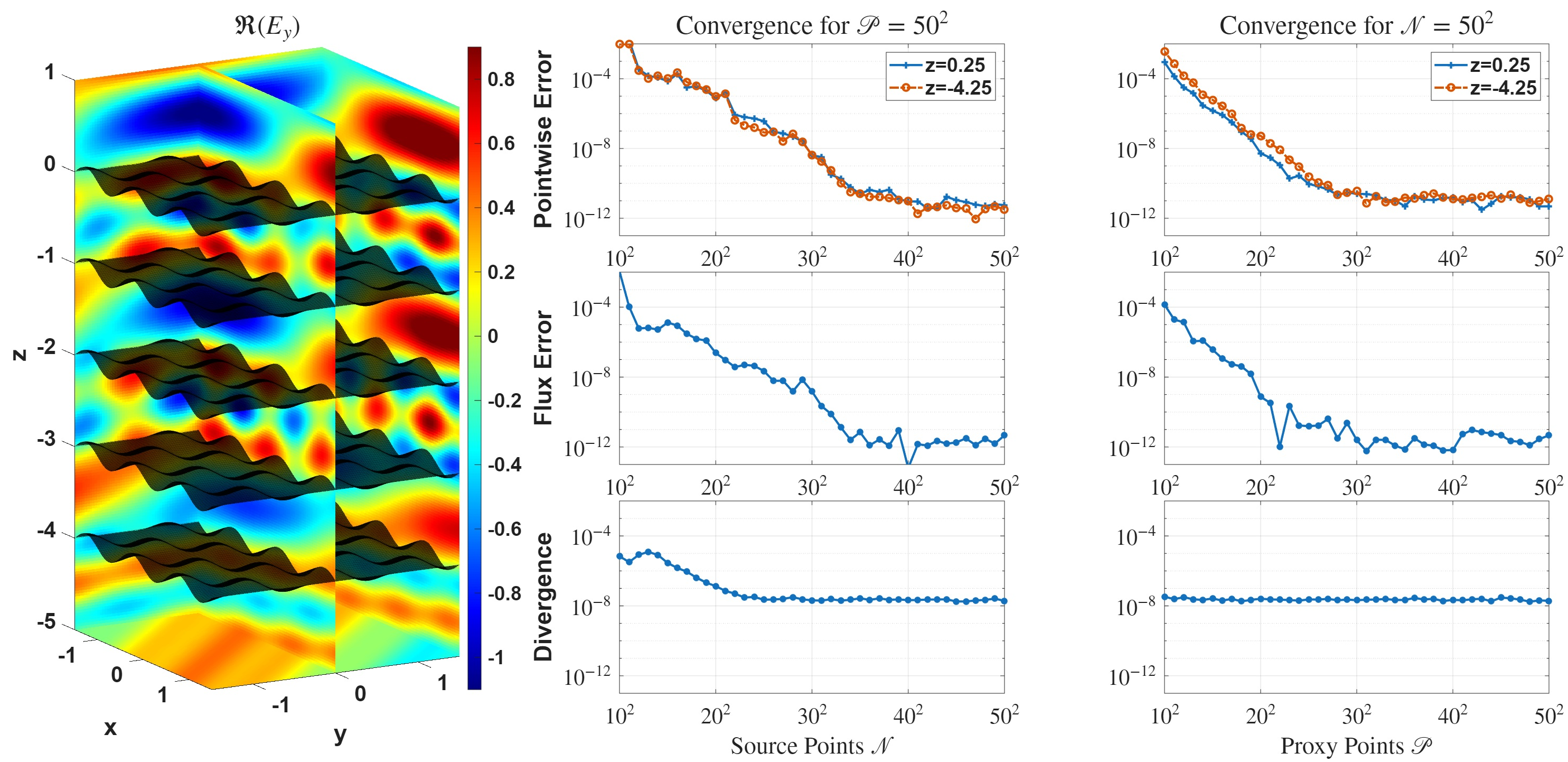} 
   \caption{Five interfaces with $f_l(x,y) = -l + 0.1\sin(2\pi x)\cos(2\pi y)$, $l=0,\ldots, 4$. \textbf{Left:} Real part of $E_y$ for the entire cell. \textbf{Center:} Pointwise convergence  at $(0,0,0.25)$ and $(0,0,-4.25)$, flux error, and divergence as a function of MFS source points with $\cP$ fixed. \textbf{Right:} Pointwise convergence  at $(0,0,0.25)$ and $(0,0,-4.25)$, flux error, and divergence as a function of proxy source points with $\cN$ fixed.}
   \label{fig:multilayerNPvary}
\end{figure}

\paragraph{Error convergence for five bi-periodic interfaces:} For the fifth numerical experiment, we define five identically-shaped surfaces as $f_l(x,y) = -l+0.1\sin(2\pi x)\cos(2\pi y)$ for $l=0,\ldots, 4$, where $\epsilon_n$ alternates between $1$ and $4$ with $\epsilon_0=1$ and $\mu_n=1$, and $\omega = 4$. The solution is evaluated at $(0,0,0.25)$ and $(0,0,-4.25)$. The TE incident field is $\EE^{\inc}(\xx) = (-\frac{\sqrt{2}}{2},\frac{\sqrt{2}}{2}, 0)e^{\iu \kk \cdot \xx}$ with a wavevector of $\kk$ pointing in the direction of azimuthal angle $\phi = 9\pi/10$ and polar angle $\theta = \pi/4$. Figure \ref{fig:multilayerNPvary} displays the real part of the electric field component $E_y$ along with the pointwise convergence in $E_y$, flux error, and divergence. $\cP=50^2$ proxy points are used for MFS source point convergence, $\cN=50^2$ source points are used for proxy point convergence, $\cR=10$ and $\cW=30^2$ are fixed across all layers and artificial surfaces. The MFS source points are offset 0.15 units above and below each surface along their surface normals.

\begin{figure}[h!] 
   \centering
   \includegraphics[width=.32\textwidth]{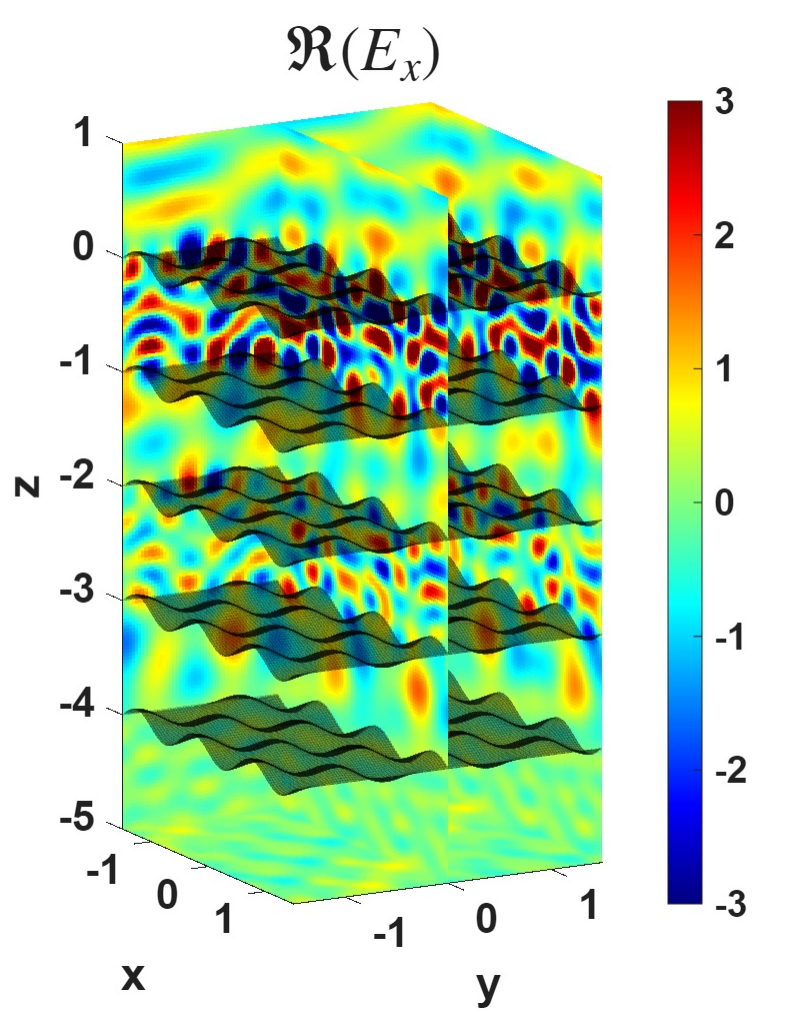}\hfill
   \includegraphics[width=.32\textwidth]{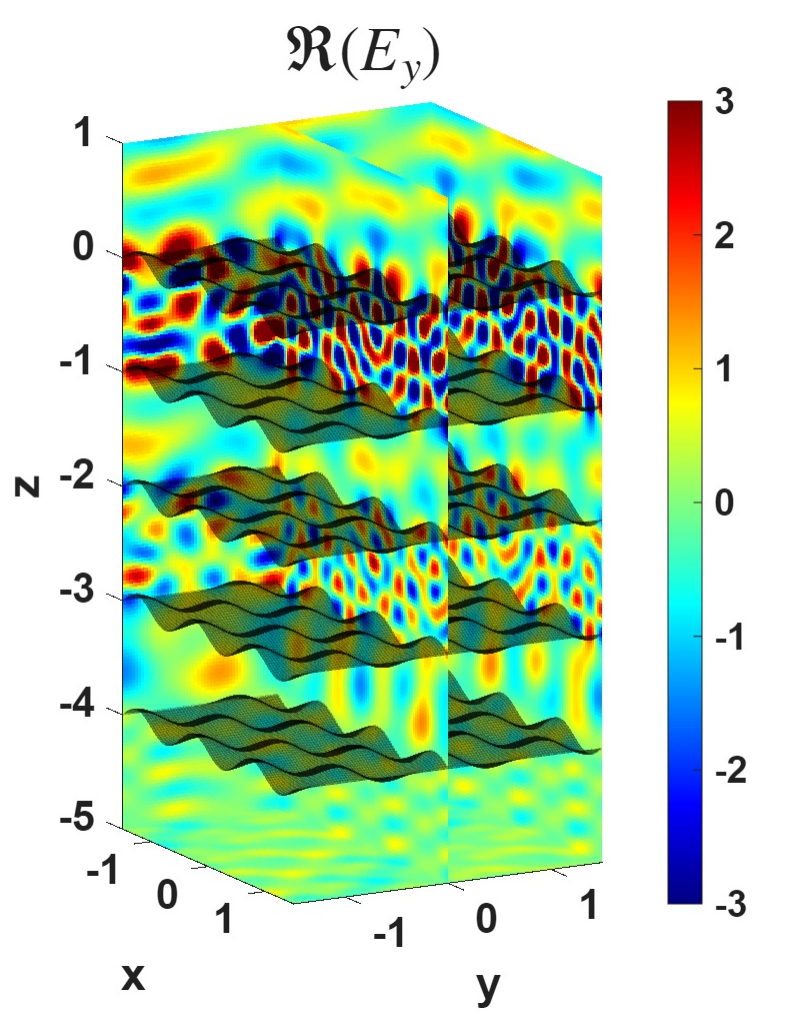}\hfill
   \includegraphics[width=.32\textwidth]{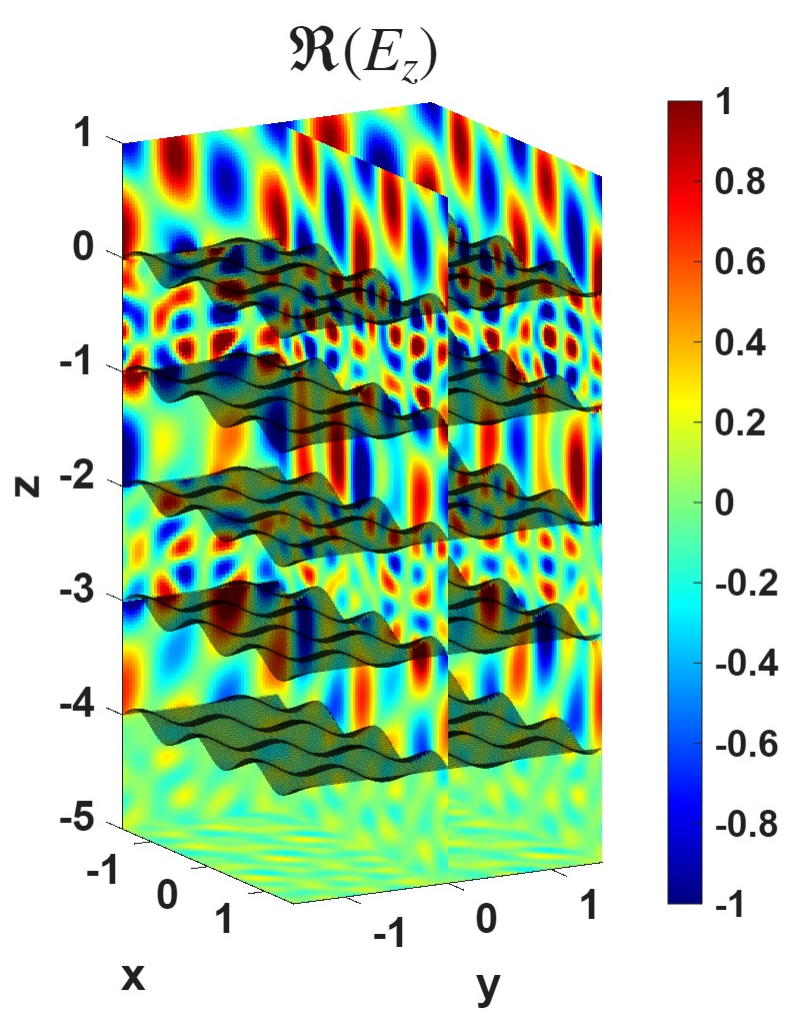}
   \caption{Five interfaces with $f_l(x,y) = -l+0.1\sin(2\pi x)\cos(2\pi y)$, $l=0,\ldots, 4$, where $\omega=10$. \textbf{Left:} Real part of $E_x$ for the entire cell. \textbf{Center:} Real part of $E_y$ for the entire cell. \textbf{Right:} Real part of $E_z$ for the entire cell.}
   \label{fig:multilayerhighomega}
\end{figure}

\paragraph{Flux error for five bi-periodic interfaces at $\omega = 10$:} For the sixth numerical experiment, we use the same settings as in the fifth numerical experiment but set $\omega = 10$. 
The real parts of the electric field $E_x$, $E_y$ and $E_z$ for the entire cell are presented in Figure \ref{fig:multilayerhighomega}. $\cN=40^2$ source points, $\cP=50^2$ proxy points, $\cR=10$ and $\cW=30^2$ are fixed across all layers and artificial surfaces. The MFS source points are offset 0.15 units above and below each respective surface along their surface normal directions, achieving a flux error of $2.19\times 10^{-10}$.

\begin{figure}[h] 
   \centering
   \includegraphics[width=.33\textwidth]{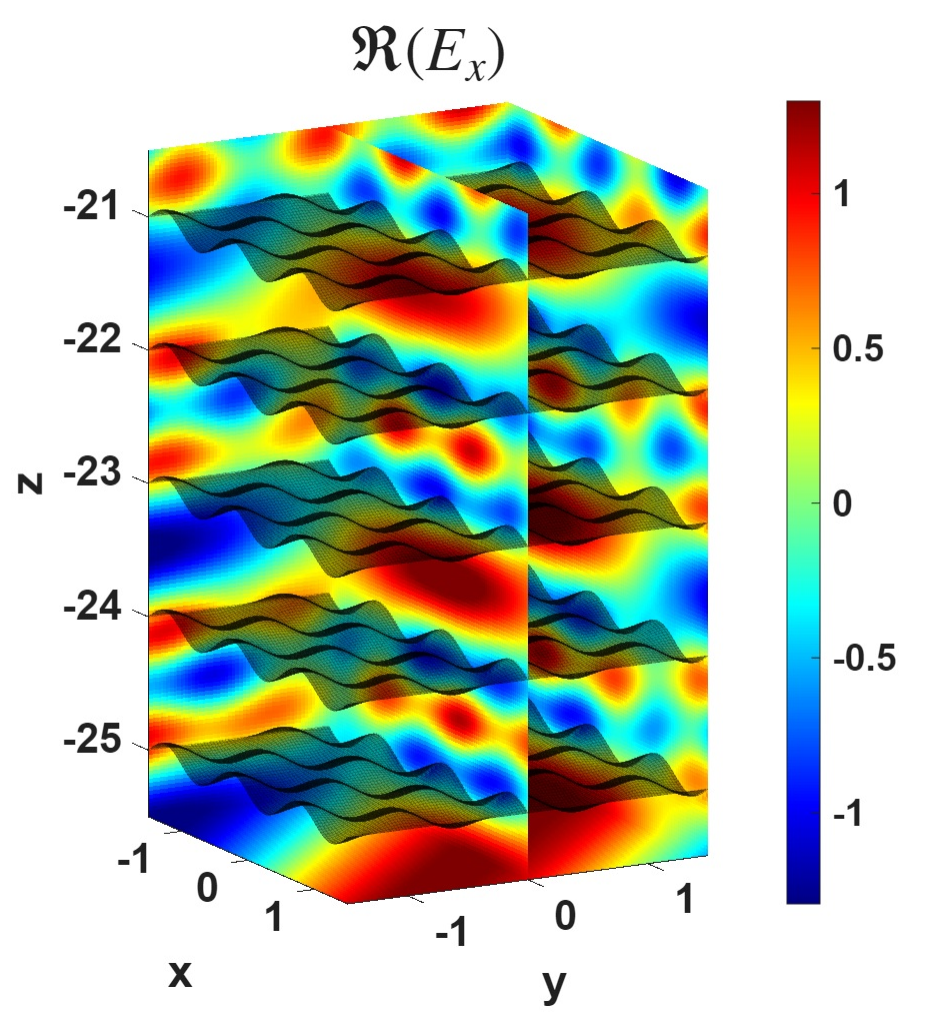}\hfill
   \includegraphics[width=.33\textwidth]{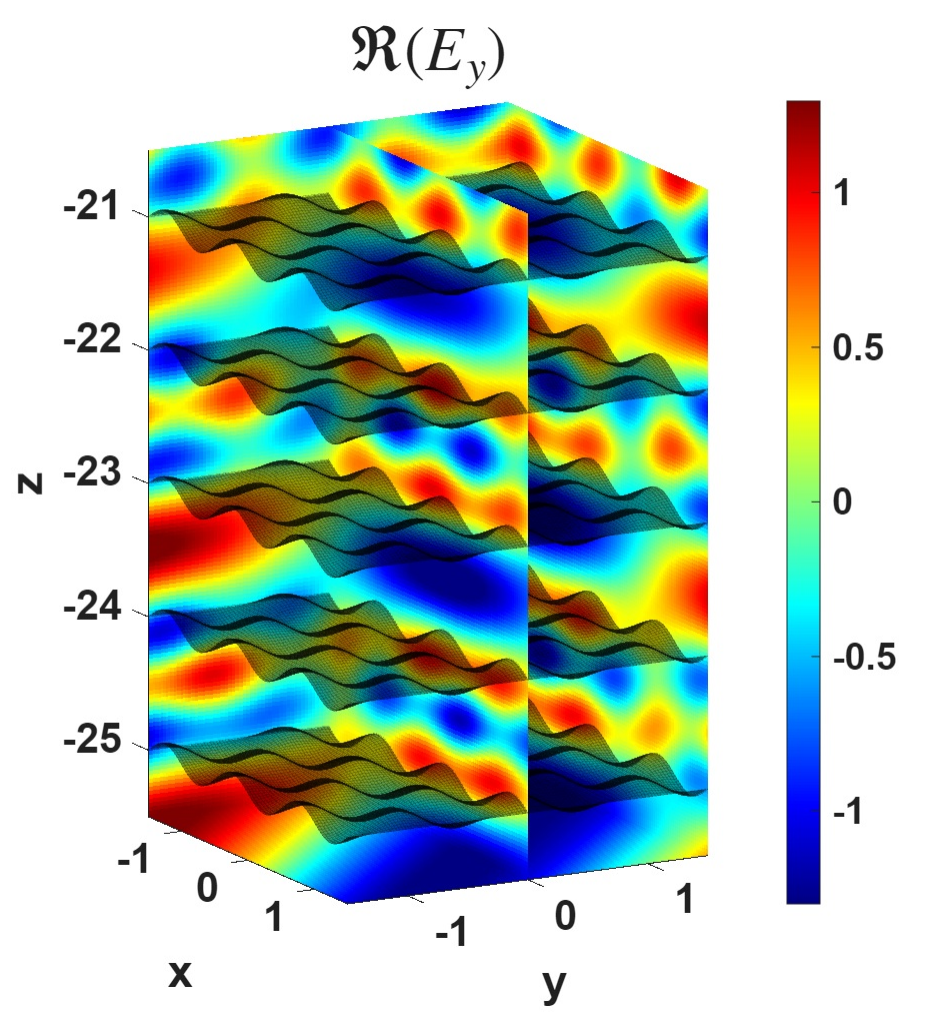}\hfill
   \includegraphics[width=.33\textwidth]{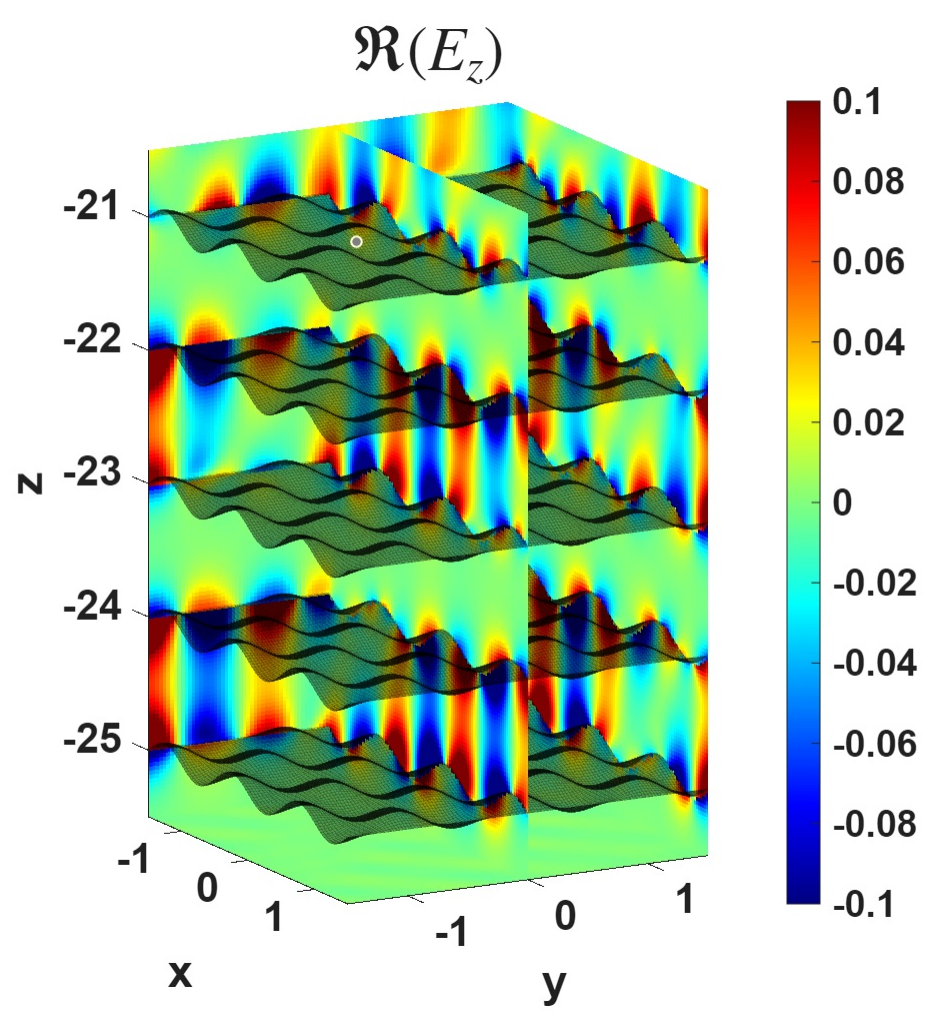}
   \caption{A portion of the field solutions between layers 22 and 27 for 40 bi-periodic interfaces defined as $f_l(x,y)=-l+0.1\sin(2\pi x)\cos(2\pi y)$, $l=0,\ldots, 39$. \textbf{Left:} Real part of the $E_x$ field solution. \textbf{Center:} Real part of the $E_y$ field solution. \textbf{Right:} Real part of the $E_z$ field solution.}
   \label{fig:multilayermanylayer}
\end{figure}

\paragraph{Flux error for many (40) bi-periodic interfaces:} For the seventh numerical experiment, we construct a cell with 40 bi-periodic interfaces defined as $f_l(x,y)=-l + 0.1\sin(2\pi x)\cos(2\pi y)$, $l=0,\ldots, 39$, where $\epsilon_n$ alternates between $1$ and $4$ with $\epsilon_0=1$ and $\mu_n=1$, and $\omega = 4$. The TE incident field $\EE^{\inc}(\xx) = (-\frac{\sqrt{2}}{2},\frac{\sqrt{2}}{2}, 0)e^{\iu \kk \cdot \xx}$ with a wavevector of $\kk$ pointing in the direction of azimuthal angle $\phi = 9\pi/10$ and polar angle $\theta = \pi/4$. Figure \ref{fig:multilayermanylayer} displays the real part of each electric field component for a portion of the entire cell. The total solution has a flux error of $1.68\times 10^{-8}$ when $\cN=30^2$ source points, $\cP=20^2$ proxy points, $\cR=10$ and $\cW=20^2$ are fixed across all layers and artificial surfaces. The MFS source points are offset 0.15 units above and below each respective surface along their surface normals.

\begin{figure}[h] 
   \centering
   \includegraphics[width=4.5in]{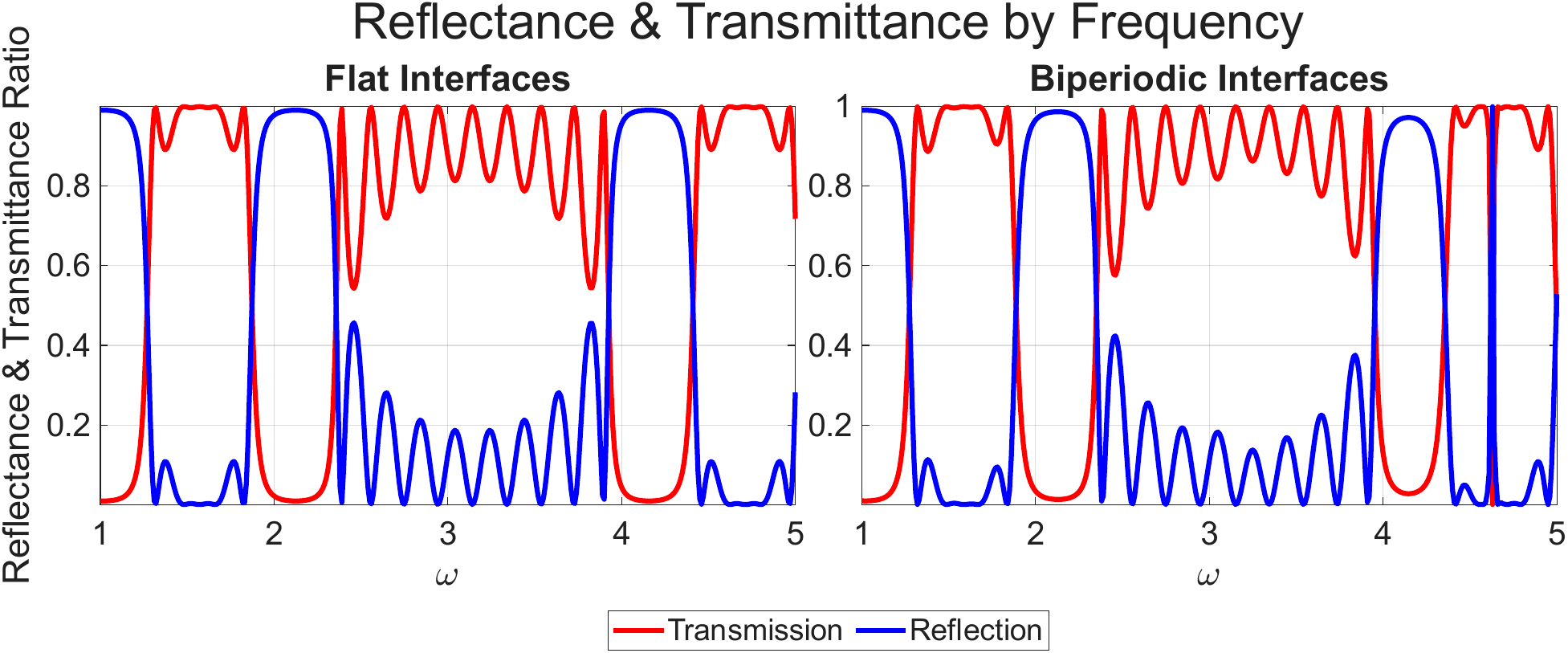} 
   \caption{Reflectance and transmittance plots for a frequency band of $\omega\in[1,5]$. \textbf{Left:} Reflectance and transmittance for a series of flat interfaces with a maximum flux error of $1.60\times10^{-8}$. \textbf{Right:} Reflectance and transmittance for a series of bi-periodic interfaces with a maximum flux error of $6.25\times10^{-5}$.}
   \label{fig:omegasweep}
\end{figure}

\paragraph{Reflectance and transmittance spectra for 10 interfaces:} For the eighth and final numerical experiment, we reuse our previous definitions of both flat and bi-periodic surfaces ($f(x,y)=0$ and $f(x,y)=0.1\sin(2\pi x)\cos(2\pi y)$) and construct cells with 10 interfaces where $\epsilon_n$ alternates between $1$ and $4$ with $\epsilon_0=1$ and $\mu_n=1$. The TE incident field $\EE^{\inc}(\xx) = (0,1, 0)e^{\iu \kk \cdot \xx}$ has a wavevector of $\kk$ pointing in the direction of azimuthal angle $\phi = \pi$ and polar angle $\theta =0$. We vary $\omega$ from $1$ to $5$ in increments of $0.01$, and the resulting transmission and reflection flux ratios are shown in Figure \ref{fig:omegasweep}. In both cases, $\cN=20^2$ source points, $\cP=20^2$ proxy points, $\cR=10$ and $\cW=20^2$ are fixed across all layers and artificial surfaces. The MFS source points are offset 0.15 units above and below each respective surface along their surface normals.

\section{Conclusion}\label{sec:conclusion}
The combination of MFS and periodization scheme yields accurate numerical solutions of the time-harmonic Maxwell's equations for multiple layers near machine precision. We demonstrate the optimal choice of MFS points, pointwise convergence with respect to the exact solution, and zero divergence conditions for a flat interface. The numerical solution agrees with the Helmholtz equation solution when the interface is periodic only in the $x$-axis. For multiple layers with bi-periodic interfaces, 
flux error convergence, pointwise convergence, and zero-divergence conditions are presented. Moreover, transmission and reflection spectra over $\omega$ of 10 interfaces are computed for use in potential applications. This work is the first step toward developing robust and fast integral equation methods for Maxwell's equations in bi-periodic multilayered media.

\section*{Acknowledgment}

J. Huang was supported by an NSF grant (DMS2152289) and B. Wu was supported by a grant  from the Simons Foundation (SFI-MPS-TSM-00013549).


%



  \bibliographystyle{elsarticle-num} 
  \bibliography{citations.bib}



%
%
%
\end{document}